\documentclass[reqno]{amsart}

\usepackage{graphicx}
\usepackage{mathrsfs}
\usepackage{mathtools}
\usepackage{amsmath}
\usepackage{caption}
\usepackage{color}
\usepackage{subcaption}
\mathtoolsset{showonlyrefs=true}
\numberwithin{equation}{section}

\newtheorem{thm}{Theorem}[section]
\newtheorem{lem}[thm]{Lemma}
\newtheorem{pro}[thm]{Proposition}

\newtheorem{cor}[thm]{Corollary}

\newtheorem{claim}{Claim}

\def\R{{\mathbb R}}
\def\Z{{\mathbb Z}}

\def\bb{\begin}
\def\be{\begin{equation}}
\def\ee{\end{equation}}
\def\bea{\begin{eqnarray}}
\def\eea{\end{eqnarray}}
\def\beaa{\begin{eqnarray*}}
\def\eeaa{\end{eqnarray*}}

\def\ifl{\iffalse}

\def\pr{{\bf pr}}

\def\bb{\begin}
           \def\ea{\end{array}}
          \def\ec{\end{center}}
     \def\ed{\end{description}}
\def\be{\bb{equation}}        \def\ee{\end{equation}}
\def\bea{\bb{eqnarray}}       \def\eea{\end{eqnarray}}
\def\beaa{\bb{eqnarray*}}     \def\eeaa{\end{eqnarray*}}
 \def\et{\end{thebibliography}}

        \def\d{\delta}
      \def\e{\varepsilon}

\def\bar2{\doublebar}

\def\qed {\hfill $\Box$\vskip5pt}

\def\Sing{{\rm Sing}}
\def\Closed{{\rm Clos}}
\def\Critical{{\rm Crit}}
\def\longrightharpoonup{\relbar\joinrel\rightharpoonup}

\def\dist{{\rm dist}}

\def\Rep{{\rm Rep}}

\def\Orb{{\rm Orb}}
\def\OS{{\rm OrientSh}}
\def\StSh{{\rm StSh}}

\def\Ps{{\rm Ps}}
\def\Pt{{\rm Pt}}

\begin{document}

\title{oriented and standard shadowing properties on closed surfaces}

\author{Sogo Murakami}



\begin{abstract}
We prove that oriented and standard shadowing properties are equivalent for topological flows on closed surfaces with the nonwandering set consisting of the finite number of critical elements (i.e., singularities or closed orbits).
Moreover, we prove that each isolated singularity of a topological flow on a closed surface with the oriented shadowing property is either asymptotically stable, backward asymptotically stable,
or admits a neighborhood which splits into two or four hyperbolic sectors.

\end{abstract}

\maketitle{}

\section{Introduction}\label{sec.intro}
In the qualitative study of dynamical systems, the shadowing property has been a key property in the hyperbolic theory (see \cite{S.P.} or \cite{R} for instance), which is useful in the proof of  the structural stability of hyperbolic systems.
One of the highlights of the studies on the relation between various shadowing properties and the structural stability or the $\Omega$-stability is due to Pilyugin and Tikhomirov, who proved that the Lipschitz shadowing property and the structural stability are equivalent for $C^1$ vector fields on a compact manifold \cite{P.S.T.}.
The Lipschitz shadowing property is stronger than the standard shadowing property,  while there are several shadowing properties weaker than the standard one.
Among them, we can find  the  oriented shadowing property.
This concerns only time-orientations of orbits, allowing  orientation preserving time-reparametrizations for true orbits in order to shadow pseudotrajectories.
It is an interesting problem to know when the oriented shadowing property implies the standard one. 

As a classical work, Komuro \cite{Komuro} proved the equivalence of oriented and standard shadowing properties for topological flows on a compact metric space without singularities. Recently, this has been extended to those with finitely many Lyapunov stable or backward Lyapunov stable singularities \cite{Mura}. When the compact metric space is changed to a closed surface, it is expected that a stronger property holds.
However, we have an important example by Tikhomirov \cite{Tikho} for which the equivalence does not hold.
The example is a  $C^1$ vector field on a four-dimensional closed manifold, whose nonwandering set consists of finitely many singularities including two hyperbolic saddles with index two.
In this paper, we prove the equivalence for topological flows on a closed surface whose nonwandering set consists of a finite number of critical elements.
Note that all hyperbolic singularities of saddle type  have index one here, so the situation that Tikhomirov constructed in his example cannot occur in this case. 

First we introduce some notation.
Let $(M, \dist)$ be a $C^0$ closed surface
with a distance function on $M$.
Take a topological flow $\phi$ on $M$ (i.e., $\phi : \R \times M \to M$ is a continuous map satisfying $\phi(0, x) = x$ and $\phi(s + t, x) = \phi(s, \phi(t, x))$ for all $s, t \in \R$ and $x \in M$).

We say that $\xi  : \R \to M$ is a {\it $d$-pseudotrajectory} of $\phi$ if
\[
\dist \bigl( \xi(t + s), \, \phi(\xi(t), s) \bigr) < d
\]
for all $t \in \R$ and $s \in [0, 1]$.
Let $\Ps(d)$ be the set of all $d$-pseudotrajectories of $\phi$.
Denote by $\Rep$ the set of all homeomorphisms from $\R$ to $\R$ which preserves the orientation.
For $\e > 0$, let
\[
\Rep(\e) = \left\{ f \in \Rep ; \left\lvert \frac{f(a) - f(b)}{a - b} - 1 \right\rvert < \e, \forall a, b \in \R, a > b \right\}.
\]
We say that a topological flow $\phi$ has the {\it standard shadowing property} if for every $\e > 0$ there exists $d > 0$ such that
if $\xi \in \Ps(d)$ then
\[
\dist \bigl( \xi(t), \, \phi(h(t), x) \bigr) < \e, \quad t \in \mathbb{R}
\]
for some $x \in M$ and $h \in \Rep(\e)$.
Denote the set of all topological flows with the standard shadowing property by $\StSh(M)$.
As a weaker form of shadowing properties, we say that a topoligocal flow $\phi$ has the {\it oriented shadowing property} if for every $\e > 0$ there exists $d > 0$ such that
if $\xi \in \Ps(d)$ then
\[
\dist \bigl( \xi(t), \, \phi(h(t), x) \bigr) < \e, \quad t \in \mathbb{R}
\]
for some $x \in M$ and $h \in \Rep$.
Denote the set of all topological flows with the oriented shadowing property by $\OS(M)$.
Singularities (i.e., fixed points for topological flows) and closed orbits of $\phi$ are called {\it critical elements}.
Let $\Sing(\phi)$ be the set of all singularities of $\phi$,
and let $\Critical(\phi)$ be that of all critical elements of $\phi$,
We say that a $\phi$-invariant compact set $K$ is {\it asymptotically stable} if for any neighborhood $V$ of $K$, there exists a neighborhood $U$ of $K$ such that
\[
\phi(t, x) \in V, \quad t \geq 0
\]
and
\[
\phi(t, x) \to p, \quad t \to \infty
\]
for all $x \in U$.
Similarly, a $\phi$-invariant compact set $K$ is said to be {\it backward asymptotically stable} if for any neighborhood $V$ of $K$, there exists a neighborhood $U$ of $K$ such that
\[
\phi(t, x) \in V, \quad t \leq 0
\]
and
\[
\phi(t, x) \to p, \quad t \to -\infty
\]
for all $x \in U$.
For $r > 0$ and $x \in M$, let $B(r, x)$ be an open ball centered at $x$ with radius $r$.
Denote the $\alpha$-limit set (resp. $\omega$-limit set) of $x$ by $\alpha(x)$ (resp. $\omega(x)$), and denote the nonwandering set of $\phi$ by $\Omega(\phi)$.

\begin{figure}
    \begin{tabular}{cc}
      \begin{minipage}[t]{0.45\hsize}
        \centering
        \includegraphics[width=3.5cm,clip]{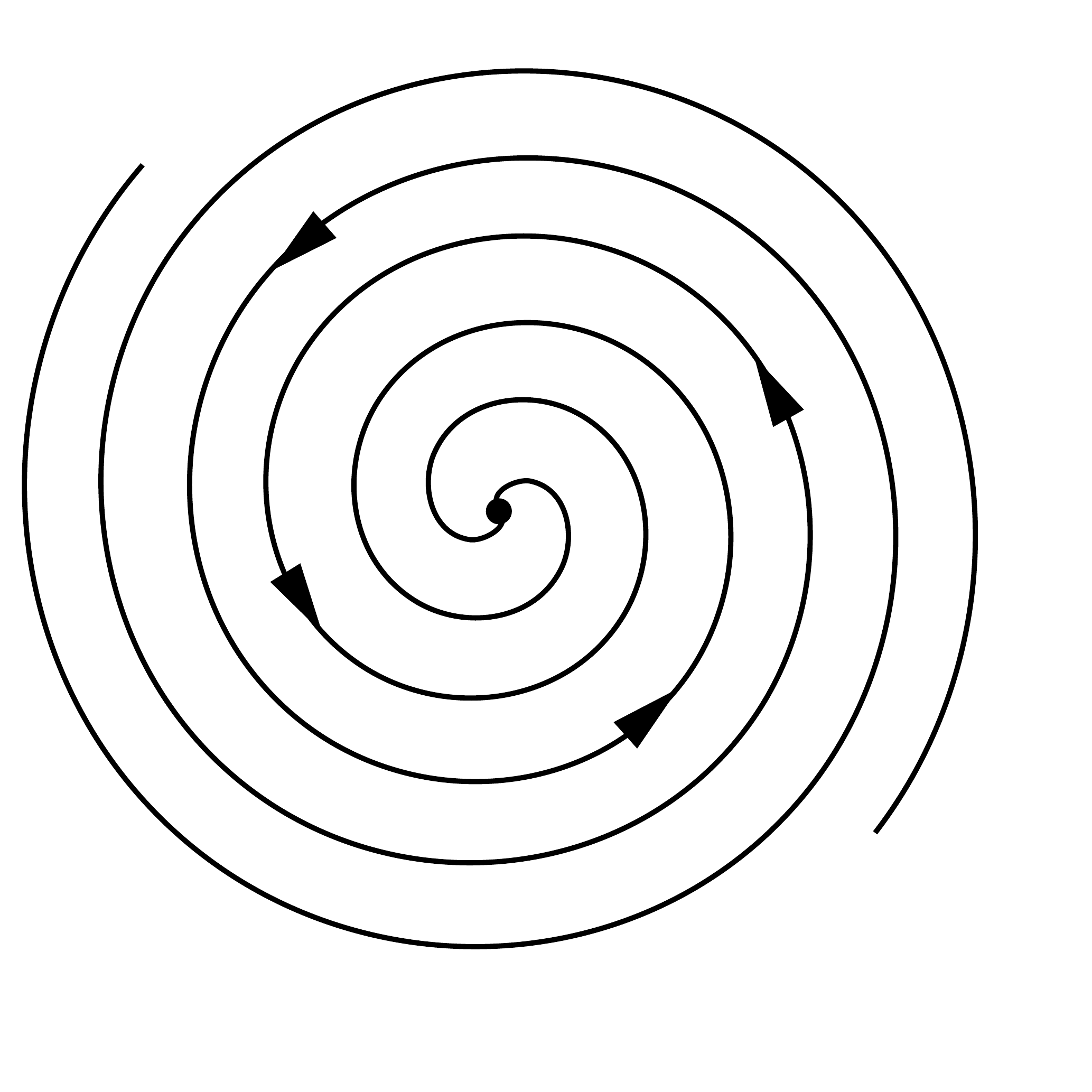}
        \subcaption{Asymptotically stable singularity.}
      \end{minipage} &
      \begin{minipage}[t]{0.45\hsize}
        \centering
        \includegraphics[width=3.5cm,clip]{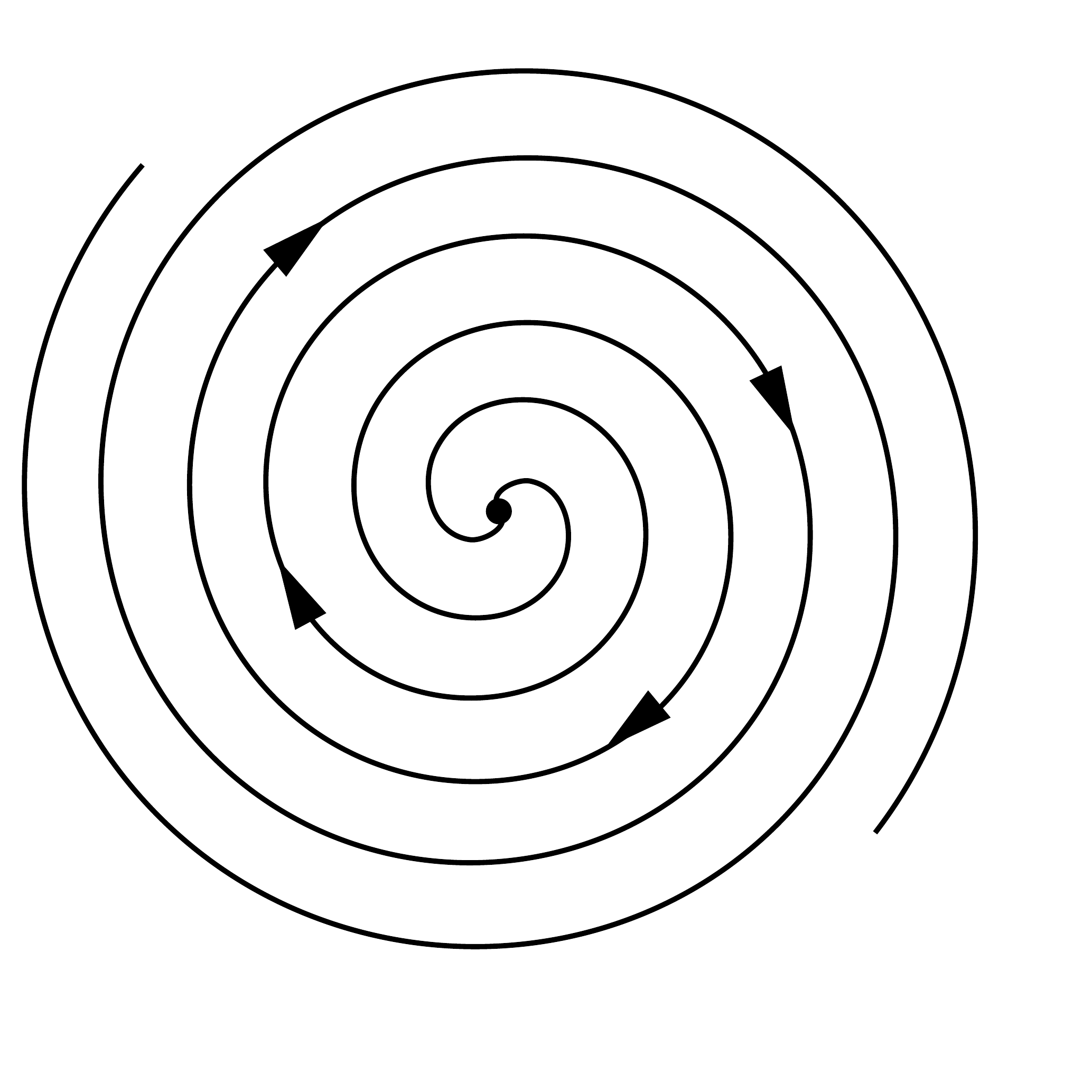}
        \subcaption{Backward asymptotically stable singularity.}
      \end{minipage}\\
      \begin{minipage}[t]{0.45\hsize}
        \centering
        \includegraphics[width=3.5cm,clip]{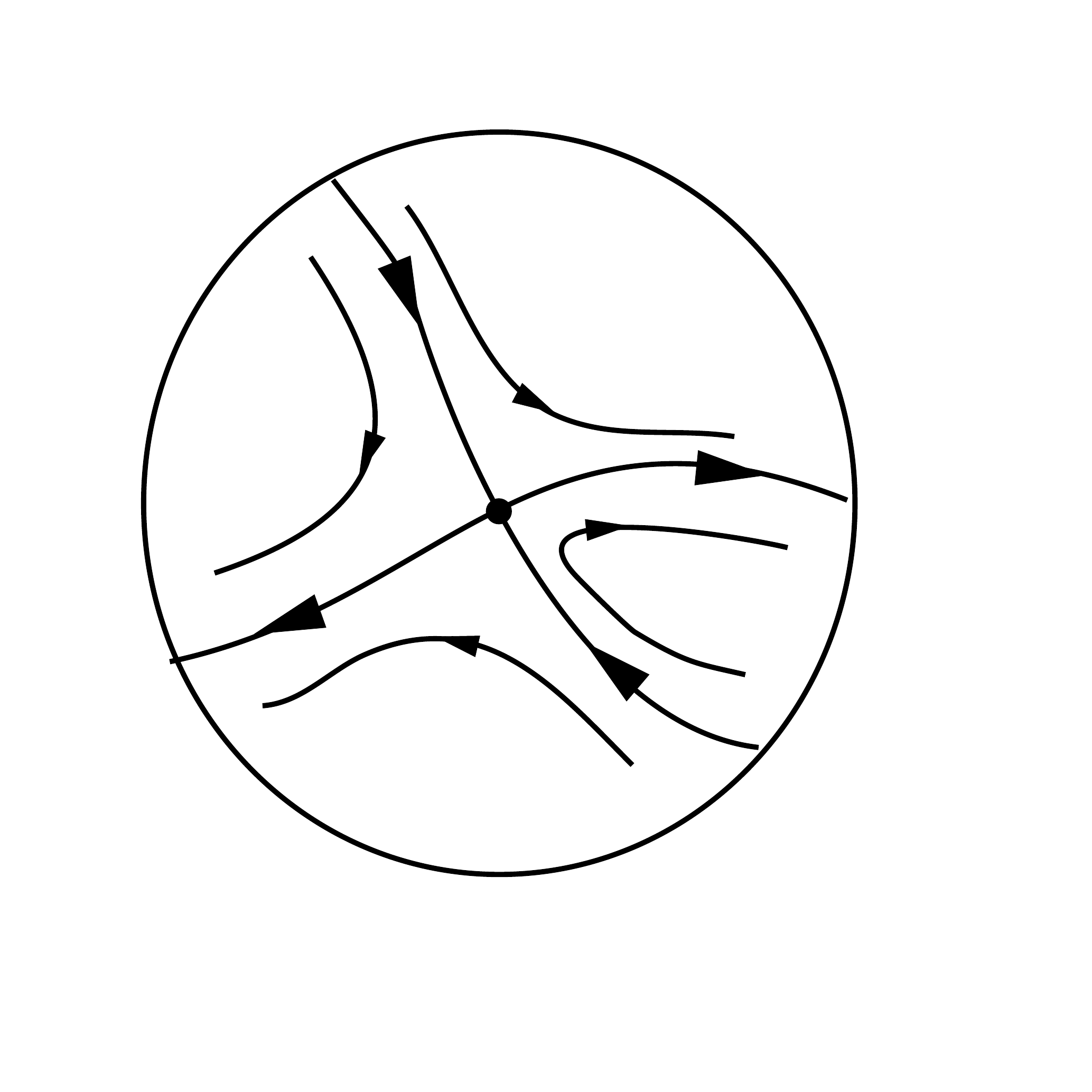}
        \subcaption{Singularity with four hyperbolic sectors.}
      \end{minipage} &
      \begin{minipage}[t]{0.45\hsize}
        \centering
        \includegraphics[width=3.5cm,clip]{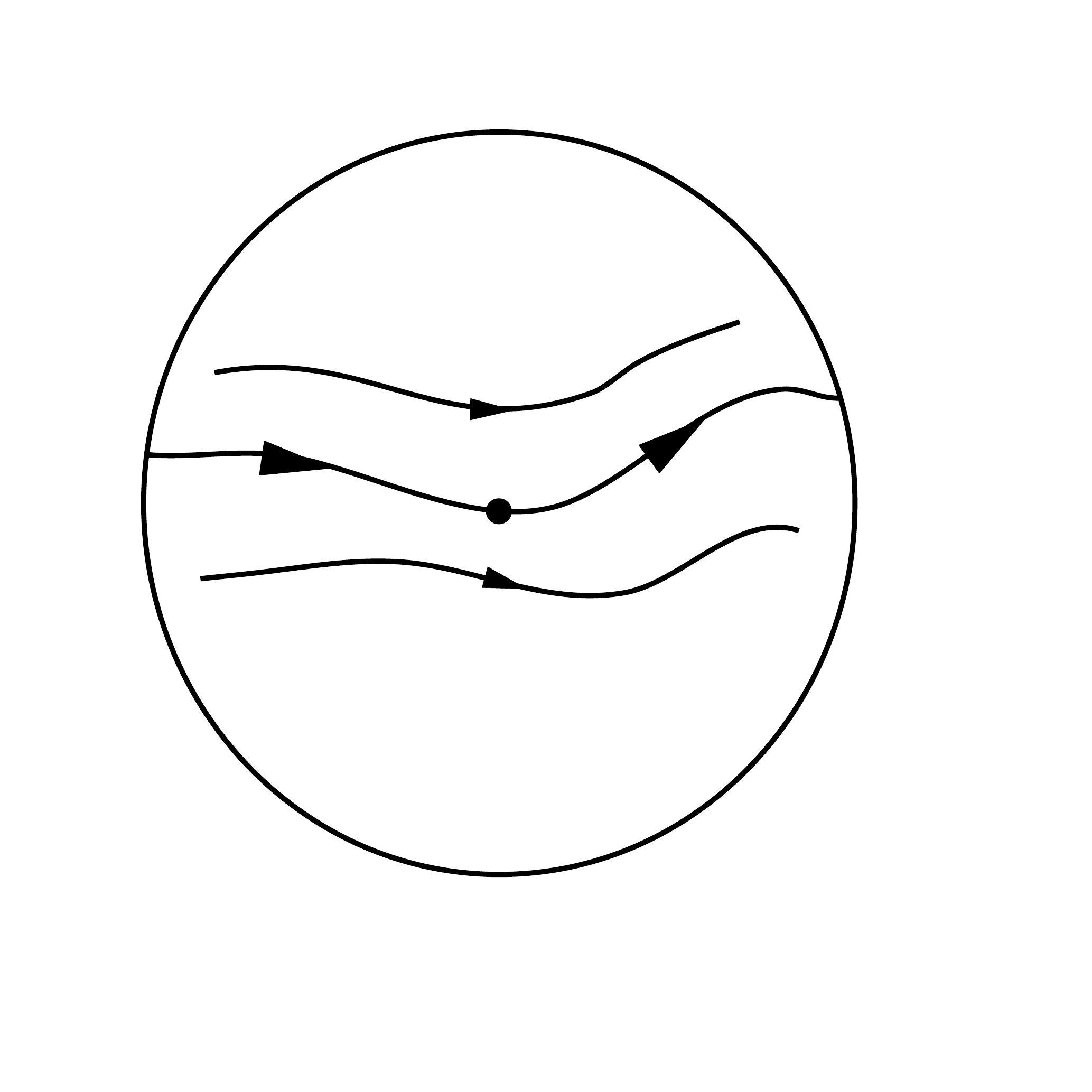}
        \subcaption{Singularity with two hyperbolic sectors.}
      \end{minipage}
    \end{tabular}
\caption{Four types of singularities in Theorem \ref{thm.classofsing}.}\label{fig.sing}
\end{figure}

Our first theorem is a classification result on singularities of $\phi \in \OS(M)$.
There are only four types of singularities, which are simple enough to understand local dynamical behavior around singularities.
A similar result can be found in \cite{MaiGu} but ours is more precise and based on the theory of sectors \cite{Hartman} that was not used in \cite{MaiGu}.
So, we provide a full proof of it for completeness.
\begin{thm}\label{thm.classofsing}
Let $M$ be a $C^0$ closed surface
and let $\phi \in \OS(M)$ with $p \in \Sing(\phi)$.
Suppose that $p$ has a neighborhood $D$ homeomorphic to a closed disk in $\R^2$ without
containing any $\alpha$-limit point or $\omega$-limit point except $p$.
Then one of the following properties hold:
\begin{itemize}
  \item[(a)] $p$ is asymptotically stable.
  \item[(b)] $p$ is backward asymptotically stable.
  \item[(c)] $D$ has an open dense subset consisting of four hyperbolic sectors associated with $\partial D$.
  \item[(d)] $D$ has an open dense subset consisting of two hyperbolic sectors associated with $\partial D$.
\end{itemize}
\end{thm}
We give definitions of hyperbolic sectors in Section \ref{sec.auxilarily} (see also Figure \ref{fig.sing}).
Using Theorem \ref{thm.classofsing}, we prove the equivalence between oriented and standard shadowing properties.
\begin{thm}\label{thm.2dim}
Let $M$ be a 
$C^0$ closed surface,
and let $\phi$ be a topological flow on $M$.
Suppose that the nonwandering set $\Omega(\phi)$ consists of a finite number of critical elements.
Then the following conditions are equivalent:
\begin{itemize}
  \item $\phi \in \StSh(M).$
  \item $\phi \in \OS(M).$
\end{itemize}
\end{thm}
It is trivial that $\StSh(M) \subset \OS(M)$.
For the proof of Theorem \ref{thm.2dim}, we need to show that $\OS(M) \subset \StSh(M)$.
The following corollary is immediate from Theorem \ref{thm.2dim}.
\begin{cor}
Let $M$ be a smooth closed surface,
and let $\phi$ be a $C^1$ Axiom A flow on $M$.
Then the following conditions are equivalent:
\begin{itemize}
  \item $\phi \in \StSh(M).$
  \item $\phi \in \OS(M).$
\end{itemize}
\end{cor}

In Section \ref{sec.auxilarily}, several propositions and lemmas are given, which play important roles in the proofs of the main theorems.
In Section \ref{sec.classofsing}, Theorem \ref{thm.classofsing} is proved.
In Section \ref{sec.prolog_thm1_2}, we reduce Theorem \ref{thm.2dim} to a proposition, and introduce some notations and lemmas.
In Sections \ref{sec.neighborhoods_orbit} and \ref{sec.neighborhoods_sing},
the behavior of pseudotrajectories near critical elements is investigated in order to prove Theorem \ref{thm.2dim}.
In Section \ref{sec.prfof2dim_2},
we finish the proof of Theorem \ref{thm.2dim}.

\section{Preliminaries}\label{sec.auxilarily}
In this section, we give several lemmas and propositions that were proved elsewhere,
which are basic to the proof of Theorems \ref{thm.classofsing} and \ref{thm.2dim}.
The first lemma is essentially proved in \cite{Bhatia}.
\begin{lem}[\cite{Bhatia}]\label{lem.pseudoflowbox}
Let $x_0 \in M \setminus \Sing(\phi)$ and let $U$ be a neighborhood of $x_0$.
Then there exist a set $S$ containing $x_0$ and $\tau_0 > 0$ such that the map
\[
\varphi : (-\tau_0, \tau_0) \times S \ni (x, t) \mapsto \phi(x, t) \in U
\]
is a homeomorphism onto its image that is an open neighborhood of $x_0$.
\end{lem}
\paragraph{\bf Remark}
The fact that ${\rm Im} \varphi$ contains a neighborhood of $x_0$ follows from an easy additional argument of the proof of \cite[Theorem 2.9]{Bhatia}, which is ommited here.

The following lemma is immediate from the Jordan-Schoenflies theorem \cite{Cairns}.
\begin{lem}[\cite{Cairns}]\label{lem.JS}
Let $C$ be a simple closed curve in $\R^2$.
Then for all $p \notin C$ and $q \in C$, there exists a continuous map $i : [0, 1] \to \R^2$ such that $i(0) = p$, $i(1) = q$ and
\[
i \bigl( [0, 1) \bigr) \cap C = \emptyset.
\]
\end{lem}
The following lemma is \cite[Theorem $31.2$, Theorem $31.5$]{Willard}.
\begin{lem}[\cite{Willard}]\label{lem.surjinj}
Let $X$ be a Hausdorff topological space,
and let $G : [-1, 1] \to X$ be a surjective continuous map with $G(-1) \neq G(1)$.
Then there exists an injective continuous map $j : [-1, 1] \to X$ such that $j(-1) = G(-1)$ and $j(1) = G(1)$.
\end{lem}
Using these lemmas, we can prove the following proposition, which claims that the flow box theorem is valid for topological flows on a closed surface.
This fact has been used in many papers (see \cite{Gut} or \cite{Aranson} for instance),
but as far as the author knows the proof has not been given explicitly.
We give a proof here for completeness.
\begin{pro}\label{pro.flowbox}
Let $x_0 \in M \setminus \Sing(\phi)$ and let $U$ be a neighborhood of $x_0$ disjoint from the limit set of $\phi$.
Then there exist an embedded closed interval $I$ with $x_0$ and $\tau_0 > 0$ such that the map
\[
\psi : (-\tau_0, \tau_0) \times I \ni (t, x) \mapsto \phi(t, x) \in U
\]
is a homeomorphism onto its image that is a neighborhood of $x_0$.
\end{pro}
\begin{prf}
By Lemma \ref{lem.pseudoflowbox},
there exist a set $S$ with $x_0$ and $\tau_0 > 0$ such that the map
\[
\varphi : (-\tau_0, \tau_0) \times S \ni (t, x) \mapsto \phi(t, x) \in M
\]
is a homeomorphism onto its image that is an open neighborhood of $x_0$.
Making $U$ smaller if necessary, we may assume that $U$ is homeomorphic to an open disk in $\R^2$ and $U \subset {\rm Im} \varphi$.
Let $C \subset U$ be a simple closed curve surrounding $x_0$.
Since $U$ and the limit set of $\phi$ are disjoint, we can take
\begin{align}
t_1 &= \sup \{ t < 0 ; \phi(t, x_0) \in C \},\\
t_2 &= \inf \{ t > 0 ; \phi(t, x_0) \in C \}.
\end{align}
Set $y_i = \phi(t_i, x_0)$ for $i = 1, 2$.
Then $C$ is divided into two subarcs $C_1$ and $C_2$ whose end points are $y_1$ and $y_2$,
and define another subarc with the same end points:
\[
C_0 = \{ \phi(t, x_0) ; t_1 < t < t_2 \}.
\]
Let $U_1$ and $U_2$ be the interiors of the simple closed curves $C_0 \sqcup C_1$ and $C_0 \sqcup C_2$, respectively.
Applying Lemma \ref{lem.JS} to $C_0 \sqcup C_1$ and $C_0 \sqcup C_2$,
we have a continuous map $i : [-1, 1] \to U$ satisfying
$i(0) = x_0$, $i(-1) \in U_1$, $i(1) \in U_2$ and 
\[
i \bigl( [-1, 0) \cup (0, 1] \bigr) \cap (C_0 \sqcup C) = \emptyset.
\]
Define $F : [-1, 1] \to S$ by
\[
F(s) = \pr_2 \circ \varphi^{-1}(i(s))
\]
for all $s \in [-1, 1]$, where $\pr_2$ is the projection to the second component.
Notice that $\pr_1 \circ \varphi^{-1}(x) \to 0$ as $x \to x_0$,
where $\pr_1$ is the projection to the first component.
For sufficiently small $\eta \in (0, 1)$,
we have
$F(-\eta) \in U_1$,
$F(\eta) \in U_2$ and
\[
F \bigl( [-\eta, 0) \cup (0, \eta] \bigr) \cap (C_0 \sqcup C) = \emptyset.
\]
Define $G : [-1, 1] \to S$ by $G(s) = F(\eta s)$.
Then apply Lemma \ref{lem.surjinj} to $G : [-1, 1] \to G([-1, 1])$ in order to have
an injective continuous map $j : [-1, 1] \to G([-1, 1])$ such that $j(-1) = G(-1)$ and $j(1) = G(1)$.
Since $G([-1, 1]) \subset U_1 \sqcup C_0 \sqcup U_2$ and $j([-1, 1])$ is connected,
there exists $t_0 \in (-1, 1)$ such that $j(t_0) \in C_0$.
Since $C_0$ is a connected subset of an orbit and $U \subset {\rm Im} \, \varphi$, we have $C_0 \cap S = \{ x_0 \}$.
Thus $j(t_0) = x_0$ is obtained.
Let $I = j([-1, 1])$ and
define $\psi : (-\tau_0, \tau_0) \times I \to {\rm Im} \, \varphi$ by
\[
\psi(t, x) = \phi(t, x).
\]
Then, by the choice of $S$, $\psi$ is injective, and the Invariance of Domain Theorem shows that $\psi$ is a homeomorphism onto its image.
Finally, $x_0$ is in the interior of ${\rm Im} \, \psi$, since $x_0$ is in an open set
$\psi( (-\tau_0, \tau_0) \times j((-1, 1) )$ of ${\rm Im} \, \psi$.
\qed
\end{prf}
Now we introduce some notations in order to give the next propsition.
Let $p$ be a singularity of $\phi$, and let $U$ be a neighborhood of $p$ homeomorphic to an open disk in $\R^2$.
Let $C \subset U$ be a simple closed curve surrounding $p$ and let $U_C$ be the interior of $C$.
For $x \in C$, we say that $\Orb^+(x)$ (resp. $\Orb^-(x)$) is a
{\it positive} $($resp. {\it negative}$)$ {\it base orbit from $x$} if
\[
\phi(t, x) \in U_C, \quad t > 0 \text{ (resp. $t < 0$) }
\]
and
\[
\phi(t, x) \to p, \quad t \to \infty \text{ (resp. $t \to -\infty$). }
\]
Let $c_i$ be a positive or negative base orbit from $x_i \in C$ for $i = 1, 2$.
Open subset $S$ of $U_C$ is called {\it sector} determined by $c_1$ and $c_2$ if the boundary of $S$ consists of $c_1, c_2$, $p$ and a subarc $C_{12}$ from $x_1$ to $x_2$.
A sector $S$ determined by $c_1$ and $c_2$ is called {\it elliptic} if
$\Orb(x_1) = \Orb(x_2) \subset \overline{S}$,
and {\it hyperbolic} if
$S$ is not elliptic and $C_{12} \cup S$ contains neither positive nor negative base orbit.
A sector $S$ determined by $c_1$ and $c_2$ is called {\it positive} $($resp. {\it negative}$)$ {\it parabolic} sector if $c_1$ and $c_2$ are positive (resp. negative) base orbits and $\overline{S}$ contains no negative (resp. positive) base orbit.
As argued in Proposition \ref{pro.nonpara} below, we may easily see that these three types of sectors do not overlap.

The following proposition is the Local Sector Classification Theorem \cite[Lemmas 8.2 and 8.3]{Hartman}.
This was proved in the $C^1$ setting, but the proof is still valid for topological flows.

\begin{figure}
\setlength{\unitlength}{20mm}
\begin{tabular}{ccc}
\begin{minipage}[t]{0.33\hsize}
\begin{picture}(0, 2)
        \put(0.1, 0){\includegraphics[width=3.5cm,clip]{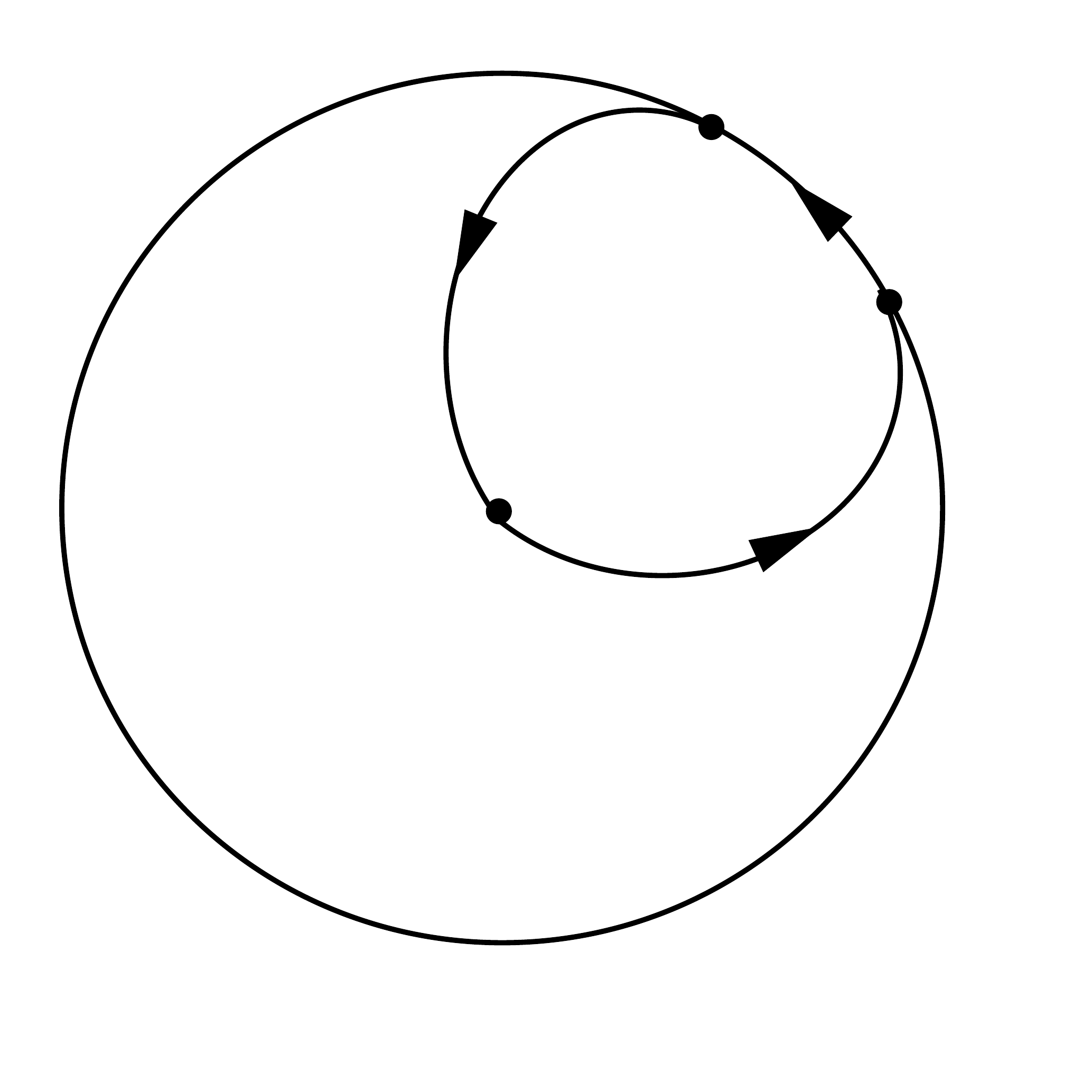}}
        \put(0.8, 0.8){\makebox(0, 0)[b]{\setlength{\fboxsep}{0.25mm}%
        	\colorbox{white}{\scriptsize $p$}}}
        \put(1.4, 1.6){\makebox(0, 0)[b]{\setlength{\fboxsep}{0.25mm}%
        	\colorbox{white}{\scriptsize $x_1$}}}
        \put(1.7, 1.3){\makebox(0, 0)[b]{\setlength{\fboxsep}{0.25mm}%
        	\colorbox{white}{\scriptsize $x_2$}}}
\end{picture}
\subcaption{Elliptic sector. \ \ \ \ }
\end{minipage}
\begin{minipage}[t]{0.33\hsize}
\begin{picture}(0, 2)
        \put(0.1, 0){\includegraphics[width=3.5cm,clip]{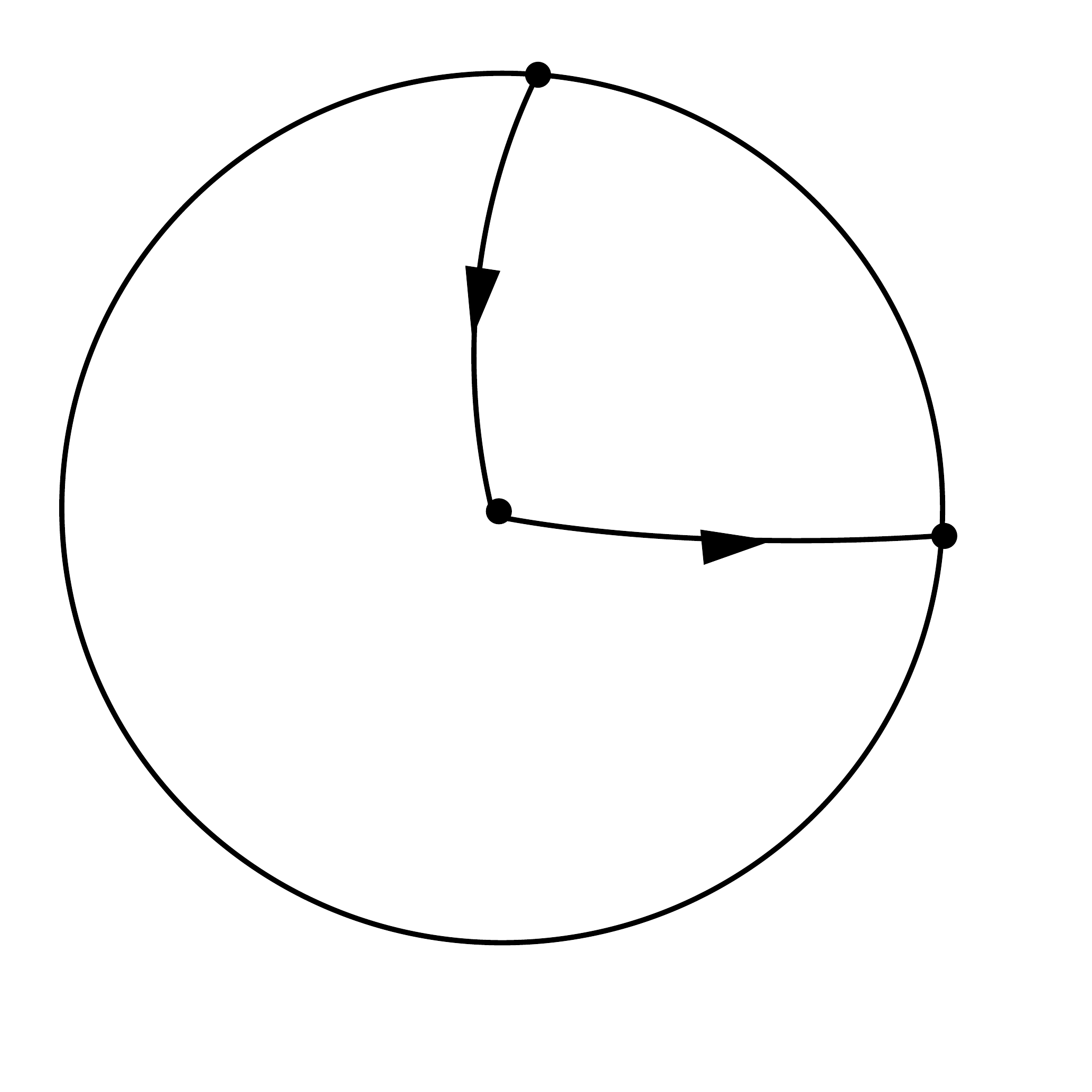}}
        \put(0.8, 0.8){\makebox(0, 0)[b]{\setlength{\fboxsep}{0.25mm}%
                \colorbox{white}{\scriptsize $p$}}}
        \put(1.1, 1.7){\makebox(0, 0)[b]{\setlength{\fboxsep}{0.25mm}%
                \colorbox{white}{\scriptsize $x_1$}}}
        \put(1.8, 0.9){\makebox(0, 0)[b]{\setlength{\fboxsep}{0.25mm}%
                \colorbox{white}{\scriptsize $x_2$}}}
\end{picture}
\subcaption{Hyperbolic sector. \ \ \ \ }
\end{minipage}
\begin{minipage}[t]{0.33\hsize}
\begin{picture}(0, 2)
        \put(0.1, 0){\includegraphics[width=3.5cm,clip]{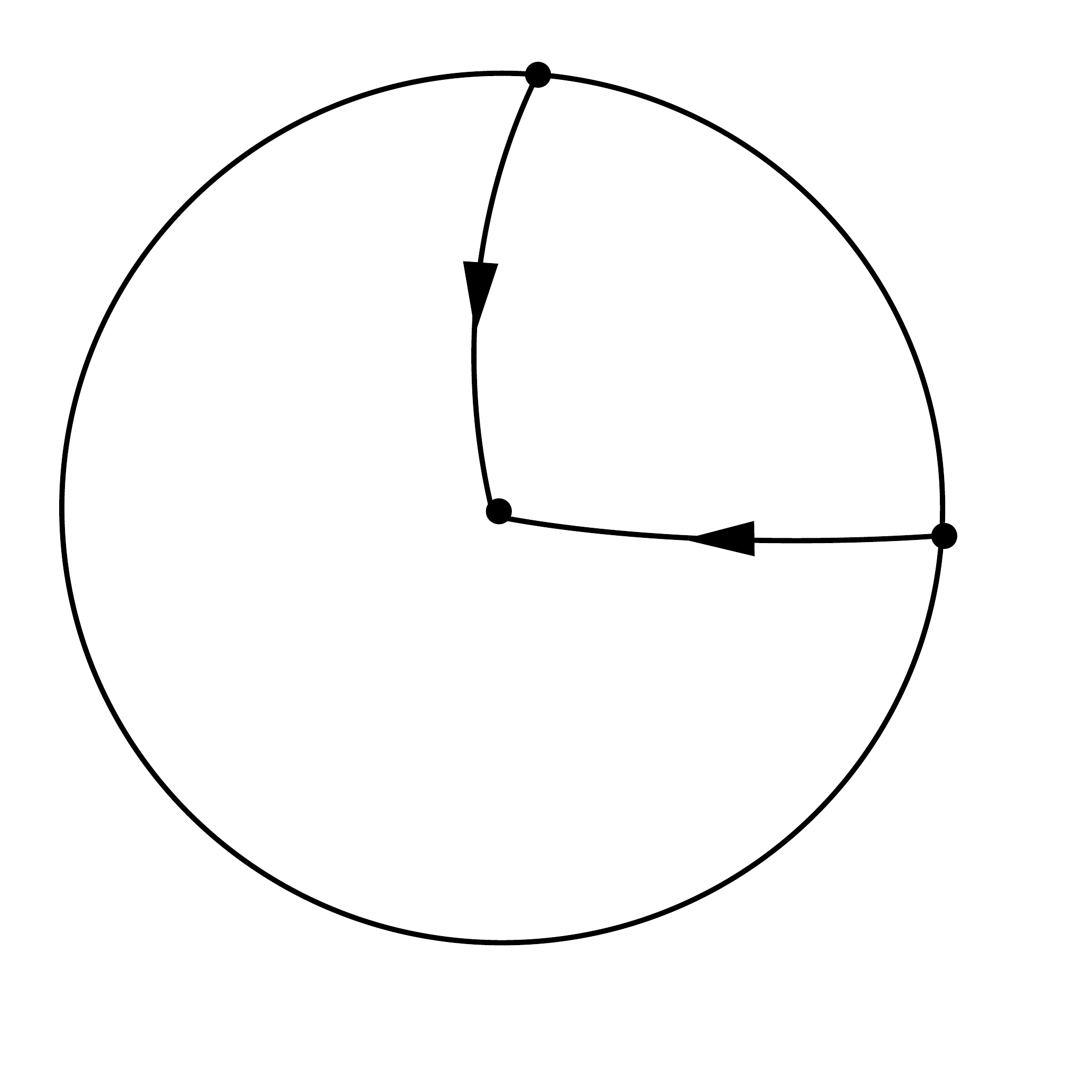}}
        \put(0.8, 0.8){\makebox(0, 0)[b]{\setlength{\fboxsep}{0.25mm}%
        	\colorbox{white}{\scriptsize $p$}}}
        \put(1.1, 1.7){\makebox(0, 0)[b]{\setlength{\fboxsep}{0.25mm}%
        	\colorbox{white}{\scriptsize $x_1$}}}
        \put(1.8, 0.9){\makebox(0, 0)[b]{\setlength{\fboxsep}{0.25mm}%
        	\colorbox{white}{\scriptsize $x_2$}}}
\end{picture}
\subcaption{Parabolic sector. \ \ \ \ }
\end{minipage}
\end{tabular}
\caption{Three types of sectors.}
\end{figure}
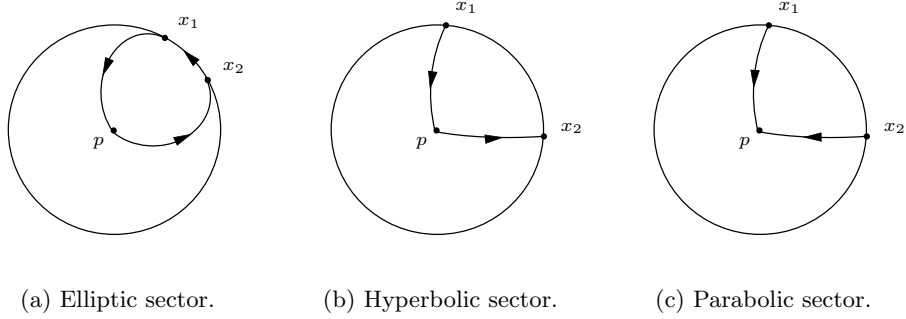
\begin{pro}[\cite{Hartman}]\label{pro.LSC}
Let $p$ be a singularity of $\phi$, and let $D$ be a closed neighborhood of $p$ homeomorphic to a closed disk in $\R^2$.
If $D$ contains no $\alpha$-limit point or $\omega$-limit point except $p$,
then $D$ has a dense subset consisting of finite number of elliptic, parabolic and hyperbolic sectors associated with $\partial D$.
\end{pro}

The three lemmas below will be used for the proof of Theorem \ref{thm.2dim}.
The following lemma has been already well-known (see \cite{R.C.} for instance).
\begin{lem}[Conley's Fundammental Theorem of Dynamical Systems]\label{lem.ordering}
There exists a continuous function $L_0 : M \to \mathbb{R}$ satisfying the following properties:
\begin{enumerate}
  \item For any $x \notin \Omega(\phi)$,
  \[
  \mathbb{R} \ni t \mapsto L_0(\phi(t, x)) \in \R
  \]
  is strictly increasing.
  \item For $x, y \in \Omega(\phi)$, $x$ and $y$ are in the same critical element iff $L_0(x) = L_0(y)$.
\end{enumerate}
\end{lem}

For $\Lambda \in \Critical(\phi)$, define
\begin{align}
W^s(\Lambda) &= \{ x \in M ; \phi(t, x) \to \Lambda, \quad t \to \infty \},\\
W^u(\Lambda) &= \{ x \in M ; \phi(t, x) \to \Lambda, \quad t \to -\infty \}.
\end{align}
For $\Lambda_1, \Lambda_2 \in \Critical(\phi)$, we write $\Lambda_1 \rightharpoonup \Lambda_2$ if
\[
W^u(\Lambda_1) \cap W^s(\Lambda_2) \setminus \Omega(\phi) \neq \emptyset.
\]
The proof of the following lemma is similar to that of  \cite[Lemma $2.4.6$]{S.P.},
which is ommitted here.
\begin{lem}\label{lem.nocycle}
If $\Omega(\phi)$ consists of finitely many critical elements,
then $\phi \in \OS(M)$ has the no-cycles condition.
\end{lem}
The following lemma is the flows version of the Birkhoff Constant Theorem \cite[Theorem 2.2.5]{S.P.}.
The proof of this lemma is similar to \cite[Theorem 2.2.5]{S.P.}.
\begin{lem}[\cite{S.P.}]\label{lem.birkhoffconst}
Let $U$ be a open set including $\Omega(\phi)$.
Then there exist $R > 0$ and $d > 0$ such that
if $x \in M$ satisfies
\[
\phi(t, x) \notin U, \quad t \in (0, \ell),
\]
then $\ell < R$.
\end{lem}

\section{Classification of singularities}\label{sec.classofsing}
In this section, we prove Theorem \ref{thm.classofsing}, which classifies topological flows on a neighborhood of each $p \in \Sing(\phi)$ with a neighborhood $D$ satisfying the condition of Theorem \ref{thm.classofsing}.

Denote by $\Sing_s(\phi)$ the set of all asymptotically stable singularities of $\phi$, and
by $\Sing_u(\phi)$ that of all asymptotically unstable ones of $\phi$.
Denote by $\Sing_4(\phi)$ the set of all singularities of $\phi$ with property $(c)$ of Theorem \ref{thm.classofsing}, and
by $\Sing_2(\phi)$ that of those with property $(d)$ of Theorem \ref{thm.classofsing}.
Take $p \in \Sing(\phi)$ and its neighborhood $D$ satisfying the condition of Theorem \ref{thm.classofsing}.
\begin{lem}\label{claim.selfconn}
For all $x \in D \setminus \{ p \}$, we have $\Orb(x) \not\subset D$.
In particular, there is no elliptic sectors associated with $\partial D$.
\end{lem}
\begin{prf}
To get a contradiction, suppose there exists a point $x_0 \in D \setminus \{ p \}$ such that $\Orb(x_0) \subset D$.
Then, by hypothesis, we have $\omega(x_0) = \alpha(x_0) = \{ p \}$.
Denote by $U$ the open region in the interior of the simple closed curve $\Orb(x_0) \cup \{ p \}$.
Since $U \subset D$,
we have
\begin{equation}
\Orb(x) \subset U \subset D\label{eq.interiororb}
\end{equation}
for all $x \in U$.
Let $x_1$ be a point in $U$ and
take $\e_1 > 0$ so small that $B(2\e_1, x_1)$ is contained in $U$.
Then \eqref{eq.interiororb} implies $\omega(x_1) = \alpha(x_1) = \{ p \}$.
For sufficiently small $d > 0$, all $\xi \in \Ps(d)$ are $\e_1$-shadowed.
There exists $T_d > 0$ such that
\[
\phi(T_d + t, x_1), \quad \phi(-T_d + t, x_1) \in B(d/2, p), \quad t \in [0, 1].
\]
Consider a $d$-pseudotrajectory $\xi : \R \to U$ such that
\[
\xi(t) = 
\phi(t - 2T_dn, x_1),\quad t \in[(2n - 1)T_d, (2n + 1) T_d)
\]
for all $n \in \Z$.
Then there exist $x_2 \in B(\e_1, x_1)$ and $h \in \Rep$ with $h(0) = 0$ such that
\[
\dist \bigl( \xi(t), \, \phi(h(t), x_2) \bigr) < \e_1, \quad t \in \R.
\]
From the choice of $x_1$ and this inequality, we have
\begin{align}
\dist \bigl( p, \, \phi(h(2T_d n), x_2) \bigr)
&\geq \dist \bigl( p, \, \xi(2T_d n) \bigr)
- \dist \bigl( \xi(2T_d n), \, \phi(h(2T_d n), x_2) \bigr)\\
&\geq 2\e_1 - \e_1 = \e_1
\end{align}
for all $n \in \Z$.
However, it follows from \eqref{eq.interiororb} that
$\omega(x_2) = \alpha(x_2) = \{ p \}$.
This is a contradiction.
\qed
\end{prf}

\begin{claim}\label{claim.notab}
If $p \notin \Sing_s(\phi)$ $($resp. $p \notin \Sing_u(\phi)$$)$, then there exists $x \in \partial D$ such that
$\Orb^-(x) \subset D$ $($resp. $\Orb^+(x) \subset D$$)$ and
\[
\phi(t, x) \to p, \quad t \to -\infty, \quad
(resp. \; \phi(t, x) \to p, \quad t \to \infty).
\]
\end{claim}
\begin{prf}
Assume $p \notin \Sing_s(\phi)$.
Then we can take $\e > 0$ such that $B(\e, p) \subset D$ and for any $n \in \Z_{> 0}$, there exists $z_n \in B(1/n, p)$ satisfying
\[
\Orb^+(z_n) \not\subset B(\e, p)
\]
or $\omega(z_n) \neq \{ p \}$.
Notice that if $\Orb^+(z_n) \subset B(\e, p)$, then $\omega(z_n) = \{ p \}$ by the hypothesis on $D$.
Therefore, for all $n \geq 1$,
\[
\Orb^+(z_n) \not\subset B(\e, p).
\]
Let
\[
t_n = \inf \{ t \geq 0 ; \phi(t, z_n) \in \partial B(\e, p) \}
\] 
and $w_n = \phi(t_n, z_n) \in \partial B(\e, p)$.
Taking a subsequence $\{ w_{n_k} \}_{k \geq 1}$ if necessary,
we may assume that $w_{n_k}$ converges to some $w \in \partial B(\e, p)$.
We have $\Orb^-(w) \subset D$ since $t_{n_k} \to \infty$ as $k \to \infty$. 
Thus by the hypothesis on $D$ again,
\[
\phi(t, w) \to p, \quad t \to -\infty.
\]
Moreover, by Lemma \ref{claim.selfconn} we have $\Orb^+(w) \not\subset D$.
Thus we may take
\[
t_0 = \inf \{ t > 0 ; \phi(t, w) \notin D \}
\]
and define $x = \phi(t_0, w)$.
Then it is easy to see that $x$ satisfies the required condition.
The proof for $p \notin \Sing_u(\phi)$ is similar reversing the time direction.
\qed
\end{prf}

\begin{pro}\label{pro.nonpara}
If $p \notin \Sing_s(\phi) \cup \Sing_u(\phi)$,
then $p$ does not have any parabolic sector associated with $\partial D$.
\end{pro}
\begin{prf}
Assume that $p$ has a positive parabolic sector $S$ determined by $\Orb^+(x_1)$ and $\Orb^+(x_2)$.
Let $C_{1, 2}$ be the subarc of $\partial D$ from $x_1$ to $x_2$ such that
$S$ is the subset of $D$ with boundary consisting of
\[
p, \Orb^+(x_1),
\Orb^+(x_2) \text{ and }
C_{1, 2}.
\]
Using Proposition \ref{pro.flowbox}, we can take an embedding $i : [-1, 1] \to D$ with $i(0) = \phi(1, x_1)$ such that an embedded closed interval $i([-1, 1])$ satisfies the condition of Proposition \ref{pro.flowbox} with $x_0 = i(0)$.
In addition, we may assume that $i((0, 1]) \subset S$.
Then
\begin{equation}
\Orb^+(i(s)) \subset S, \quad s \in (0, s_0] \label{eq.triangle}
\end{equation}
for some small $s_0 \in (0, 1]$.
In fact, suppose by contradiction that there exists a sequence $\{ s_n \}_{n \geq 1}$ in $(0, 1]$ such that
$s_n \to 0$ as $n \to \infty$ and
\[
\Orb^+(i(s_n)) \not\subset S.
\]
Then for all $n \geq 1$, there exists $t_n > 0$ such that
\[
t_n = \inf \{ t > 0 ; \phi(t, i(s_n)) \in C_{1, 2} \}.
\]
Taking a subsequence if necessary, we may assume $\phi(t_n, i(s_n)) \to y_0$ for some $y_0 \in C_{1, 2}$.
Since $t_n \to \infty$, by continuity, $\phi(t, y_0) \in \overline{S}$ for all $t < 0$,
contradicting that $S$ is parabolic.
Thus \eqref{eq.triangle} holds.
Making $s_0$ smaller if necessary, we may assume that
$\Orb^+(i(s_0)) \cap i([0, s_0]) = \{ i(s_0) \}$.
Denote by $S_0$ the region in the interior of the closed curve consisting of
\[
i([0, s_0]),
p,
\Orb^+(i(0)) \text{ and }
\Orb^+(i(s_0)),
\]
where note that $\omega(i(s_0)) = \{ p \}$ by \eqref{eq.triangle}.
Let $x_0 = \phi(1, i(s_0/2)) \in S_0$.
Take $\e_1 > 0$ satisfying the following properties:
\begin{equation}
B(\e_1, x_0) \subset S_0\label{eq.nonparapro1}
\end{equation}
and
\begin{equation}
\dist(\overline{S_0}, \partial D) > \e_1.\label{eq.nonparapro2}
\end{equation}
Since $p \notin \Sing_s(\phi)$, using Claim \ref{claim.notab},
we have $w \in \partial D$ with $\Orb^-(w) \subset D$ and
\[
\phi(-t, w) \to p, \quad t \to \infty.
\]
Our hypthesis $\phi \in \OS(M)$ implies that there exists $d > 0$ such that
every $d$-pseudotrajectory is $\e_1$-oriented shadowed.
Take $T > 0$ so large that the function $\xi : \R \to M$ defined by
\[
\xi(t) = 
\begin{cases}
\phi(t, x_0), & t \leq T,\\
\phi(t - 2T, w), & t > T
\end{cases}
\]
is a $d$-pseudotrajectory.
Then there exist $h \in \Rep$ with $h(0) = 0$ and $z_0 \in S$ such that
\[
\dist \bigl( \xi(t), \, \phi(h(t), z_0) \bigr) < \e_1, \quad t \in \R.
\]
From this inequality, we have $\dist(x_0, z_0) < \e_1$. 
Thus, by \eqref{eq.triangle} and \eqref{eq.nonparapro1}, $\phi(h(t), z_0) \in S_0$ for all $t \geq 0$.
We also have
\[
\dist \bigl( w, \, \phi(h(2T), z_0) \bigr) < \e_1.
\]
This contradicts \eqref{eq.nonparapro2} and proves  this proposition when $p$ has a positive parabolic sector.
The proof for the case where $p$ has a negative parabolic sector is similar.
\qed
\end{prf}
Assume that $p \notin \Sing_s(\phi) \cup \Sing_u(\phi)$.
By Proposition \ref{pro.LSC},
there is an open dense subset of $D$ consisting of finite number of elliptic, hyperbolic and parabolic sectors.
Then, combining Lemma \ref{claim.selfconn} and Proposition \ref{pro.nonpara},
we can observe that $D$ has an open dense subset consisting of even number of hyperbolic sectors.
\begin{lem}\label{lem.hypsecnbd}
Let $S$ be a hyperbolic sector of $p$ determined by $\Orb^+(x_1)$ and $\Orb^-(x_2)$ for some $x_1, x_2 \in \partial D$.
Let $x_1' \in \Orb^+(x_1) \cup \{ p \} \setminus \{ x_1 \}$ and $x_2' \in \Orb^-(x_2)$.
For any neighborhood $U_2$ of $x_2'$, there exists a neighborhood $U_1$ of $x_1'$ such that for any $x \in U_1 \cap S$,
\begin{equation}
\phi(t, x) \in S, \quad t \in [0, t_x]\label{eq.hypsecnbd1}
\end{equation}
and
\begin{equation}
\phi(t_x, x) \in U_2\label{eq.hypsecnbd2}
\end{equation}
for some $t_x > 0$.
\end{lem}
\begin{prf}
Since $S$ is a hyperbolic sector,
for all $x \in S$,
\begin{equation}
\Orb^+(x) \not\subset S.\label{eq.hypsecnbd3}
\end{equation}
In fact, $S$ has no positive base orbit, so if \eqref{eq.hypsecnbd3} is false 
then we would have $\Orb^-(x) \subset S$ (note that $\Orb^+(x) \subset S$ shows that $\omega(x) = \{ p \}$).
This contradicts Lemma \ref{claim.selfconn}.

Suppose to the contrary that there exists a sequence $\{ y_n \}_{n \geq 1}$ in $S$ 
such that
\[
y_n \to x_1', \quad n \to \infty
\]
and for all $n \geq 1$ and $s > 0$, either
\begin{equation}
\phi([0, s], y_n) \not\subset S\label{eq.hypsecnbd4}
\end{equation}
or
\begin{equation}
\phi(s, y_n) \notin U_2\label{eq.hypsecnbd5}
\end{equation}
hold.
By \eqref{eq.hypsecnbd3}, there exists
$t_n = \inf \{ t > 0 ; \phi(t, y_n) \in \partial S \} < \infty$.
If $\phi(s, y_n) \in U_2$ for some $s \in [0, t_n)$,
then \eqref{eq.hypsecnbd4} and \eqref{eq.hypsecnbd5} do not hold for $s$.
Thus we get
\begin{equation}
\phi([0, t_n], y_n) \cap U_2 = \emptyset \label{eq.disU}
\end{equation}
for all $n \geq 1$.
Taking a subsequence if necessary, we may assume that $\phi(t_n, y_n) \to y_0$ for some $y_0 \in \partial D$,
where $\phi(t, y_0) \in \overline{S}$ for all $t \leq 0$
since $t_n \to \infty$.
Combine this and the fact that $S$ is a hyperbolic sector to have $y_0 = x_2$.
This contradicts \eqref{eq.disU}.
\qed
\end{prf}
\begin{pro}\label{pro.prohyp}
If $p \notin \Sing_s(\phi) \cup \Sing_u(\phi)$,
then the number of hyperbolic sectors is either $2$ or $4$.
\end{pro}
\begin{prf}
Suppose to the contrary that there exist more than four hyperbolic sectors.
Let $S_1$ and $S_2$ be adjacent sectors associated with $\partial D$,
where $S_i$ is determined by $\Orb^+(x_0)$ and $\Orb^-(y_i)$, $i = 1, 2$ for some $x_0, y_1, y_2 \in \partial D$.
Since there are more than four hyperbolic sectors,
we have $y_3 \in \partial D \setminus \{ y_1, y_2 \}$ with $\Orb^-(y_3) \subset D$.
Let $\e_1 > 0$ be such that
\begin{equation}
\e_1 < \dist \bigl( \{ \phi(-1, y_1), \phi(-1, y_2) \}, \, \Orb^+(x_0) \cup \{ p \} \cup \Orb^-(y_3) \bigr)\label{eq.defe_1}
\end{equation}
and
\begin{equation}
\e_1 < \dist \bigl( y_3, \overline{S_1} \cup \overline{S_2} \bigr).\label{eq.defe_12}
\end{equation}
Using Lemma \ref{lem.hypsecnbd} with repect to $S_1$ and $S_2$, we can choose $\e_2 \in (0, \e_1)$ so small that if
$x \in B(\e_2, \phi(1, x_0)) \setminus \Orb^+(x_0)$,
then there exists $t_0 > 0$ with
\[
\phi(t, x) \in S_1 \cup S_2, \quad t \in [0, t_0]
\]
and
\[
\phi(t_0, x) \in B(\e_1, y_1) \cup B(\e_1, y_2).
\]
Let $d > 0$ be such that every $\xi \in \Ps(d)$ is $\e_2$-oriented shadowed.
The function $\xi : \R \to M$ defined by
\[
\xi(t) =
\begin{cases}
\phi(1 + t, x_0), & t \leq T_d,\\
\phi(t - 2T_d, y_3), & t > T_d
\end{cases}
\]
is a $d$-pseudotrajectory for some $T_d > 0$ and hence $\xi$ is $\e_2$-shadowed.
Take $h \in \Rep$ with $h(0) = 0$ and $z_0 \in M$ satisfying
\begin{equation}
\dist \bigl( \xi(t), \, \phi(h(t), z_0) \bigr) < \e_2, \quad t \in \R.\label{eq.prohypdist}
\end{equation}
Then we have $z_0 \in B(\e_2, \phi(1, x_0))$.

Now let us consider two possible cases.
\vspace{3mm}
\paragraph{\bf Case 1:} $z_0 \in B(\e_2, \phi(1, x_0)) \setminus \Orb^+(x_0)$.
\vspace{3mm}

By the choice of $\e_2$, there exists $t_0 > 0$ such that
\[
\phi(t, z_0) \in S_1 \cup S_2, \quad t \in [0, t_0]
\]
and
$\phi(t_0, z_0) \in B(\e_1, y_1) \cup B(\e_1, y_2)$.
Let $t_1 = \min \{ h^{-1}(t_0), 2T_d \}$.
Then
\[
\phi(h(t_1), z_0) \in \overline{S_1} \cup \overline{S_2}
\]
and
\[
\xi(t_1) \in \Orb^+(x_0) \cup \{ p \} \cup \Orb^-(y_3).
\]
In addition, we have either
\[
\xi(t_1) = y_3 \text{ when } t_1 = 2T_d
\]
or
\[
\phi(h(t_1), z_0) \in B(\e_1, y_1) \cup B(\e_1, y_2) \text{ when } t_1 = h^{-1}(t_0).
\]
By these facts together with \eqref{eq.defe_1} and \eqref{eq.defe_12}, 
\[
\dist \bigl( \xi(t_1), \, \phi(h(t_1), z_0) \bigr) > \e_1.
\]
This contradicts \eqref{eq.prohypdist} with $\e_2 < \e_1$.
\vspace{3mm}
\paragraph{\bf Case 2:} $z_0 \in B(\e_2, \phi(1, x_0)) \cap \Orb^+(x_0)$.
\vspace{3mm}

Since $\phi(h(2T_d), z_0) \in \overline{S_1} \cup \overline{S_2}$, \eqref{eq.defe_12} implies
\[
\dist \bigl( y_3, \, \phi(h(2T_d), z_0) \bigr) > \e_1.
\]
This contradicts \eqref{eq.prohypdist} with $\e_2 < \e_1$ again.
\qed
\end{prf}
Now combining Proposition \ref{pro.nonpara} with Proposition \ref{pro.prohyp}, we obtain Theorem \ref{thm.classofsing}.

\section{Reduction of Theorem \ref{thm.2dim}}\label{sec.prolog_thm1_2}
In this section, we reduce Theorem \ref{thm.2dim} to a proposition which is easier to show, and introduce some notations and lemmas.
These notations and lemmas plays an important role in the rest of this paper.
Assume that $\phi \in \OS(M)$ has the nonwandering set consisting of finite number of critical elements.

In order to prove Theorem \ref{thm.2dim}, it is enough to show that for every $\e_0 > 0$, there exists $\d > 0$ such that for all $\d$-pseudotrajectory $\xi$ we have
\begin{equation}
\dist \bigl( \xi(t), \, \phi(h(t), x) \bigr) < 2\e_0, \quad t \in \mathbb{R}\label{eq.thm2dimfirst}
\end{equation}
for some $x \in M$ and $h \in \Rep(2\e_0)$.
We may assume that $\e_0 \in (0, 4/5)$ satisfies
\begin{equation}
3\e_0 < \min \{ \dist(\Lambda_1, \Lambda_2) ; \Lambda_1, \Lambda_2 \in \Critical(\phi), \Lambda_1 \neq \Lambda_2 \}.\label{eq.defe_0}
\end{equation}
This choice of $\e_0$ will be used in the following sections.
Set
\[
K = M \setminus \bigcup_{x \in \Sing(\phi)} B(\e_0/2, x).
\]
Then we can take $T_0 > 0$ such that 
\[
\phi(t, x) \neq x, \quad t \in [0, 2T_0]
\]
for all $x \in K$ because $K \cap \Sing(\phi) = \emptyset$ and $K$ is compact.
Define
\[
\widetilde{K} = \bigcap_{t \in [-2T_0, 2T_0]} \phi(t, K).
\]
Making $T_0 > 0$ smaller if necessary, we may assume that
\begin{equation}
B(\e_0, \Lambda) \subset \widetilde{K}\label{eq.hyptildeK}
\end{equation}
for every closed orbit $\Lambda$.
Denote by $\Pt_{T_0}$ the set of all functions $\xi : \R \to M$ satisfying
\[
\xi(t + nT_0) = 
\phi(t, \xi(nT_0))
\]
for all $t \in [0, T_0)$ and $n \in \Z$. 
For $d > 0$, let
\[
\Pt_{T_0}(d) = \{ \xi \in \Pt_{T_0} ; \dist \bigl( \phi(T_0, \xi(nT_0) ), \, \xi((n + 1)T_0) \bigr) < d, \quad n \in \Z \}.
\]
The proof of Theorem \ref{thm.2dim} is redused to showing the following proposition:
\begin{pro}\label{pro.finalpurpose}
For $\e_0 > 0$ given by \eqref{eq.defe_0}, there exists $d > 0$ such that if $\xi \in \Pt_{T_0}(d)$ then
\begin{equation}
\dist \bigl( \xi(t), \, \phi(h(t), x) \bigr) < \e_0, \quad t \in \mathbb{R}\label{eq.finalpurpose}
\end{equation}
for some $x \in M$ and $h \in \Rep(\e_0)$.
\end{pro}
Let us show how Theorem \ref{thm.2dim} follows from Proposition \ref{pro.finalpurpose}.
It is easy to see that
there exists $\d > 0$ such that if $\xi$ is a $\d$-pseudotrajectory then
\[
\dist \bigl( \xi(t), \, \xi'(t) \bigr) < \e_0, \quad t \in \R
\]
for some $\xi' \in \Pt_{T_0}(d)$.
Then, by Proposition \ref{pro.finalpurpose},
we have
\[
\dist \bigl( \xi'(t), \, \phi(h(t), x) \bigr) < \e_0, \quad t \in \mathbb{R}
\]
for some $x \in M$ and $h \in \Rep(\e_0)$.
Therefore
\begin{align}
\dist \bigl( \xi(t), \, \phi(h(t), x) \bigr)
&\leq \dist \bigl( \xi(t), \, \xi'(t) \bigr)
+ \dist \bigl( \xi'(t), \, \phi(h(t), x) \bigr)\\
&< \e_0 + \e_0 = 2\e_0, \quad t \in \mathbb{R},
\end{align}
proving \eqref{pro.finalpurpose} from which Theorem \ref{thm.2dim} follows.

The two lemmas below will be needed in the proof of Proposition \ref{pro.finalpurpose}, which involve the above choice of $\e_0$.
The following lemma is essentially \cite[Proposition 2.1]{Mura}
where $\widetilde{K}$ was defined by given $T_0 \in (0, 1)$ and a compact set $K$,
while we here fixed $T_0$ depending on $\e_0$ given above.
In what follows, we write $\Pt_{T_0} = \Pt$ and
$\Pt_{T_0}(d) = \Pt(d)$.
\begin{lem}[\cite{Mura}]\label{lem.shadowing}
For every $\e > 0$, there exists $\e' > 0$ satisfying the following property:

If $y \in M$ satisfies
\[
\dist \bigl( \xi(t + t_0), \, \phi(g(t), y) \bigr) < \e', \quad t \in [0, T_1]
\]
and
\[
\xi(t + t_0), \, \phi(g(t), y) \in \widetilde{K}, \quad t \in [0, T_1]
\]
for some $\xi \in \Pt$, $t_0 \in \R$, $T_1 \geq T_0$ and $g \in \Rep$ with $g(0) = 0$,
then there exists $\tilde{g} \in \Rep(\e)$ with $\tilde{g}(0) = g(0)$ and $\tilde{g}(T_1) = g(T_1)$ such that
\[
\dist \bigl( \xi(t + t_0), \, \phi(\tilde{g}(t), y) \bigr) < \e, \quad t \in [0, T_1].
\]
\end{lem}
The following lemma is the closed orbits version of \cite[Lemma 5.1]{Mura}.
\begin{lem}[\cite{Mura}]\label{lem.genstable}
Let $\Lambda$ be an asymptotically stable $($resp. backward asymptotically stable$)$ critical element of $\phi$.
Then there exist $d > 0$ and a neighborhood $U$ of $\Lambda$ such that
if $\xi \in \Pt(d)$ satisfies $\xi(t_0) \in U$ for some $t_0 \in \R$,
then
\[
\xi(t) \in B(\e_0/2, \Lambda), \quad t \geq t_0 \: (\text{resp. $t \leq t_0$}).
\]
\end{lem}
The proof of this lemma is similar to that of \cite[Lemma 5.1]{Mura},
which is ommited here.
In fact, the proof of \cite[Lemma 5.1]{Mura} is also valid for closed orbits just by replacing singularities by them.

\section{Neighborhoods of Closed Orbits}\label{sec.neighborhoods_orbit}
In this section, we construct neighborhoods of closed orbits in order to investigate the behavior of pseudotrajectories.
\begin{lem}\label{lem.classclosedtraj}
Let $\Lambda$ be a closed orbit of $\phi$.
Then $\Lambda$ is asymptotically stable or backward asymptotically stable.
\end{lem}
\begin{prf}
Let $\e_0$ be the constant given by \eqref{eq.defe_0}.
Let $\varphi : [-1, 1] \to B(\e_0, \Lambda)$ be an embedding such that an embedded closed interval $\varphi([-1, 1])$
satisfies the condition of Proposition \ref{pro.flowbox} with $x_0 = \varphi(0)$.
For some $\d_0 > 0$, we can take a Poincar\'e map $P : \varphi([-\d_0, \d_0]) \to \varphi([-1, 1])$.
Define $\widetilde{P} : [-\d_0, \d_0] \to [-1, 1]$ by $\widetilde{P} = \varphi^{-1} \circ P \circ \varphi$.
It is enough to show that $0$ is attracting or repelling for $\widetilde{P}$.
Notice that $P$ is a homeomorphism onto its image, and thus $\widetilde{P}$ is increasing or decreasing.

Let us consider the following two cases separately.
  \vspace{2mm}\\
  {\bf Case 1:} $\widetilde{P}$ is increasing.

  It follows from \eqref{eq.defe_0} that $\widetilde{P}(t) \neq t$ hold for all $t \neq 0$.
  Since $\widetilde{P}$ is continuous, we may assume that $\widetilde{P}(a) > a$ for all $a \in (0, \d_0)$.
  The other case where $\widetilde{P}(a) < a$ for all $a \in (0, \d_0)$ is proved similarly.
  If $\widetilde{P}(a) < a$ for all $a \in (-\d_0, 0)$, then $0$ is repelling for $\widetilde{P}$.

  To prove by contradiction, let us assume that $\widetilde{P}(a) > a$ for some $a \in (-\d_0, 0)$ which implies $\widetilde{P}(a) > a$ for all $(-\delta_0, 0)$ because of \eqref{eq.defe_0}.
  Let $p_0 = \varphi(-\d_0/2)$ and let $p_1 = \varphi(\d_0/2)$.
  Take $\e > 0$ so small that the closure of $B(\e, p_0)$ is disjoint from $\Lambda$.
  Then let $\Gamma$ be the closure of $\Orb^+(B(\e, p_0))$.
  Making $\e$ smaller if neccesary, we may assume that
  \begin{equation}
  B(\e, p_1) \cap \Gamma = \emptyset.\label{eq.p_1ass}
  \end{equation}
  Since $\phi \in \OS(M)$, there exists $d > 0$ such that every $d$-pseudotrajectory is $\e$-oriented shadowed.
  Take $\tau_0, \tau_1 > 0$ so large that the function $\xi : \R \to M$ defined by
  \[
  \xi(t) =
  \begin{cases}
  \phi(t + \tau_0, p_0), & t \leq 0,\\
  \phi(t - \tau_1, p_1), & t > 0
  \end{cases}
  \]
  satisfies $\xi \in \Ps(d)$.
  By the choice of $d$,
  $\xi$ is $\e$-shadowed.
  Thus there exist $h \in \Rep$ with $h(-\tau_0) = 0$ and $x_0 \in M$ such that
  \begin{equation}
  \dist \bigl( \xi(t), \, \phi(h(t), x_0) \bigr) < \e, \quad t \in \R.\label{eq.distclosedtraj}
  \end{equation}
  Then $\dist \bigl( p_0, \, x_0 \bigr) < \e$, which together with the definition of $\Gamma$ implies that
  \[
  \phi(h(t), x_0) \in \Gamma
  \]
  for all $t \geq -\tau_0$.
  In particular, $\phi(h(\tau_1), x_0) \in \Gamma$ and then \eqref{eq.p_1ass}
  contradicts \eqref{eq.distclosedtraj} with $t = \tau_1$.
  \vspace{2mm}\\
  {\bf Case 2:} $\widetilde{P}$ is decreasing.

  Take $\d_1 \in (0, \d_0)$ satisfying $\widetilde{P}((-\d_1, \d_1)) \subset (-\d_0, \d_0)$.
  Then it follows from \eqref{eq.defe_0} that $\widetilde{P}^2(t) \neq t$ for all $t \neq 0$.
  Now assume that $\widetilde{P}^2(a) > a$ for all $a \in (0, \d_1)$.
  The proof for the other case where $\widetilde{P}(a) < a$ for all $a \in (0, \d_0)$ is similar.
  Since $\widetilde{P}$ is continuous and $\widetilde{P}^2(\d_1/2) > \d_1/2$, there exists $a_0 \in (\widetilde{P}(\d_1/2), 0)$ such that $\widetilde{P}(a_0) > \d_1/2$.
  Then $(a_0, \d_1/2)$ satisfies $\widetilde{P}^{-1}(a_0, \d_1/2) \subset (a_0, \d_1/2)$,
  thus by \eqref{eq.defe_0}, $(a_0, \d_1/2)$ is a repelling neighborhood of $0$ for $\widetilde{P}$.
\qed
\end{prf}
Denote by $\Closed_s(\phi)$ the set of all asymptotically stable closed orbits,
and denote by $\Closed_u(\phi)$ that of all backward asymptotically stable ones.
\begin{claim}\label{claim.closedclosed}
Let $\Lambda$ be a closed orbit of $\phi$ and let $\e_0$ be the constant given by \eqref{eq.defe_0}.
There exists $r_0 \in (0, \e_0/2)$ such that if $x, y \in \Lambda$ satisfy $\dist(x, y) < r_0$,
then
\[
\dist \bigl( \phi(t, x), \, \phi(t, y) \bigr) < \frac{\e_0}{8}, \quad t \in \R.
\]
\end{claim}
\begin{prf}
Choose $\e > 0$ so small that
\begin{equation}
\dist \bigl( \phi(t, x), \, x \bigr) < \frac{\e_0}{8}\label{eq.deffore}
\end{equation}
for all $t \in (-\e, \e)$ and $x \in \Lambda$.
To prove by contradiction, suppose that there exist sequences $\{ x_n \}_{n \geq 1}$ and $\{ y_n \}_{n \geq 1}$ in $\Lambda$ such that
\begin{equation}
\dist \bigl( x_n, y_n \bigr) < \frac{1}{n}\label{eq.x_ndist}
\end{equation}
and
\[
\dist \bigl( \phi(t_n, x_n), \, \phi(t_n, y_n) \bigr) \geq \frac{\e_0}{8}
\]
for some $\{ t_n \}_{n \geq 1}$.
Then by \eqref{eq.deffore},
\[
\phi(t_n, y_n) \not\in \phi((-\e, \e), \phi(t_n, x_n)),
\]
implying
\begin{equation}
y_n \notin \phi \bigl( (-\e, \e), \, x_n \bigr)\label{eq.nbdofx_n}
\end{equation}
for $n \geq 1$.
Taking a subsequence if necessary, we may assume that $x_n \to x_0$ for some $x_0 \in \Lambda$.
Then, by \eqref{eq.nbdofx_n}, $y_n \notin \phi \bigl( (-\e/2, \e/2), x_0 \bigr)$ for sufficiently large $n$.
This contradicts \eqref{eq.x_ndist} for large $n$.
\qed
\end{prf}
The following proposition is the main result of this section, whose proof is given after the proof of Lemma \ref{lem.orbitversion} below.
The proposition shows that when $\Lambda$ is asymptotically stable it has a neighborhood in which
not only any forward pseudotrajectory with sufficiently small jump stays in the neighborhood,
but also it is $\e_0$-standard shadowed by the forward orbit of every nearby point.
\begin{pro}\label{pro.closed}
Let $\Lambda \in \Closed_s(\phi)$ $($resp. $\Lambda \in \Closed_u(\phi)$$)$.
Then there exist a neighborhood $U_\Lambda \subset B(\e_0, \Lambda)$ of $\Lambda$ and $d_\Lambda, r_\Lambda > 0$ such that if $\xi \in \Pt(d_\Lambda)$ satisfies $\xi(t_0) \in U_\Lambda$ for some $t_0 \in \R$ then the following properties hold:
\begin{enumerate}
  \item[(a)]
  $\xi(t) \in B(\e_0/2, \Lambda)$ for all $t \geq t_0$ $($resp. $t \leq t_0$$)$.
  \item[(b)]
  For all $x \in B(r_\Lambda, \xi(t_0))$, there exists $h \in \Rep(\e_0)$ with $h(0) = 0$ such that
  \[
  \dist \bigl( \xi(t), \, \phi(h(t - t_0), x) \bigr) < \e_0, \quad t \geq t_0
   \: (\text{resp. $t \leq t_0$$)$}.
  \]
\end{enumerate}
\end{pro}
In order to prove Proposition \ref{pro.closed}, we need the following lemma.
This lemma claims that every forward orbit near a asymptotically stable closed orbit is standard shadowed by a sufficiently close forward orbit.
\begin{lem}\label{lem.orbitversion}
Let $\Lambda \in \Closed_s(\phi)$ $($resp. $\Lambda \in \Closed_u(\phi)$$)$.
Then there exists $r > 0$ such that if $x_0 \in \Lambda$ and $x_1 \in B(r, x_0)$,
then
\[
\dist \bigl( \phi(t, x_0), \, \phi(g(t), x_1) \bigr) \leq \frac{\e_0}{4}, \quad t \geq 0 \: (\text{resp. $t \leq 0$})
\]
for some $g \in \Rep(\e_0/4)$ with $h(0) = 0$.
\end{lem}
\begin{prf}
Without loss of generality, we may assume $\Lambda \in \Closed_s(\phi)$,
since the case for $\Lambda \in \Closed_u(\phi)$ is similar.
Let $r_0 \in (0, \e_0/2)$ be the constant given by Claim \ref{claim.closedclosed}. 
Apply Lemma \ref{lem.shadowing} with $\e$ replaced by $\min \{ \e_0/8, r_0/2 \}$ and let $\e_0' \in (0, r_0)$ be the constant given by the proposition for $\e = \min \{ \e_0/8, r_0/2 \}$.
We may take  $r \in (0, r_0/2)$ such that every $\xi \in \Pt(r)$ is $\e_0'$-oriented shadowed and
\[
\Orb^+(x) \subset B(\e_0/8, \Lambda)
\]
for all $x \in B(r, \Lambda)$ since $\Lambda$ is asymptotically stable.
For every $x_0 \in \Lambda$ and $x_1 \in B(r, x_0)$, set
\[
\xi(t) =
\begin{cases}
\phi(t, x_0), & t < 0,\\
\phi(t, x_1), & t \geq 0.
\end{cases}
\]
Since $\xi \in \Pt(r)$,
there exist $h \in \Rep$ with  $h(0) = 0$ and $x_1' \in M$ such that
\begin{equation}
\dist \bigl( \xi(t), \, \phi(h(t), x_1') \bigr) < \e_0', \quad t \in \R.
\end{equation}
By the choice of $r$,
\begin{equation}
\xi(t) \in B(\e_0/8, \Lambda), \quad t \in \R. \label{eq.bydefofxi}
\end{equation}
This fact and \eqref{eq.hyptildeK} enable us to apply Lemma \ref{lem.shadowing}.
Thus we have
\begin{equation}
\dist \bigl( \xi(t), \, \phi(g(t), x_1') \bigr) < \min \{ \e_0/8, r_0/2 \}, \quad t \in \R
\label{eq.orientfor}
\end{equation}
for some $g \in \Rep(\min \{ \e_0/8, r_0/2 \})$ with  $g(0) = 0$.
This inequality and \eqref{eq.bydefofxi} imply that
\[
\Orb(x_1') \subset B(\e_0/4, \Lambda).
\]
Thus, by \eqref{eq.defe_0}, we have $\alpha(x_1') = \omega(x_1') = \Lambda$.
Then Lemma \ref{lem.nocycle} shows $x_1' \in \Lambda$.
Since \eqref{eq.orientfor} implies
\[
\dist \bigl( x_1', \, x_0 \bigr)
\leq \dist \bigl( x_1', \, x_1 \bigr)
+ \dist \bigl( x_1, \, x_0 \bigr)
< r_0/2 + r < r_0,
\]
we have, $x_1' \in \Lambda \cap B(r_0, x_0)$.
Therefore, by Claim \ref{claim.closedclosed} and \eqref{eq.orientfor},
\begin{align}
\dist \bigl( \phi(t, x_0), \, \phi(g^{-1}(t), x_1) \bigr)
&\leq \dist \bigl( \phi(t, x_0), \, \phi(t, x_1') \bigr)\\
&+ \dist \bigl( \phi(t, x_1'), \, \phi(g^{-1}(t), x_1) \bigr)\\
&= \frac{\e_0}{8}
+ \dist \bigl( \phi(g(g^{-1}(t)), x_1'), \, \xi(g^{-1}(t)) \bigr)\\
&\leq \frac{\e_0}{8} + \frac{\e_0}{8} = \frac{\e_0}{4}
\end{align}
for all $t \geq 0$.
Combining this and $g^{-1} \in \Rep(\e_0/4)$, we finish the proof.
\qed
\end{prf}
Now, let us prove Proposition \ref{pro.closed}.
\vspace{3mm}\\
{\bf Proof of Proposition \ref{pro.closed}}.
Let us prove the proposition when $\Lambda \in \Closed_s(\phi)$, since the case where $\Lambda \in \Closed_u(\phi)$ can be proved similarly.
Let $r_0, r > 0$ be the positive numbers given by Claim \ref{claim.closedclosed}
and Lemma \ref{lem.orbitversion}, respectively.
Apply Lemma \ref{lem.shadowing} with $\e$ replaced by $\min \{ \e_0/4, r/2 \}$ and let $\e' \in (0, r_0)$ be the constant given by the proposition for $\e = \min \{ \e_0/4, r/2 \}$.
Let $d_\Lambda$ be such that every $\xi \in \Pt(d_\Lambda)$ is $\e'$-oriented shadowed.
Set $r_\Lambda = r/2$ and
\begin{equation}
U_\Lambda = B(r/2, \Lambda).\label{eq.defUL}
\end{equation}
Lemma \ref{lem.genstable} implies that, by making $d_\Lambda$ and $U_\Lambda$ smaller if necessary,
we may assume that if $\xi \in \Pt(d_\Lambda)$ satisfies $\xi(t_0) \in U_\Lambda$ for some $t_0 \in \R$, then
\begin{equation}
\xi(t) \in B(\e_0/2, \Lambda)\label{eq.prfclosed}
\end{equation}
for all $t \geq t_0$, proving $(a)$.
As for the proof of $(b)$, as in the hypothesis, suppose that $\xi \in \Pt(d_\Lambda)$ satisfies $\xi(t_0) \in U_\Lambda$ for some $t_0 \in \R$.
Let $x_1 \in B(r_\Lambda, \xi(t_0))$.
By \eqref{eq.defUL}, we can choose $x_0 \in \Lambda$ satisfying $\xi(t_0) \in B(r/2, x_0)$.
On the other hand, by the choice of $d_\Lambda$, there exist $x_2$ and $h \in \Rep$ with $h(0) = 0$ such that
\[
\dist \bigl( \xi(t + t_0), \, \phi(h(t), x_2) \bigr) < \e', \quad t \in \R.
\]
Since $\e' < \e_0/2$, together with \eqref{eq.prfclosed} and \eqref{eq.hyptildeK} we may use Lemma \ref{lem.shadowing} again.
Then we have
\begin{equation}
\dist \bigl( \xi(t + t_0), \, \phi(g(t), x_2) \bigr) < \min \{ \e_0/4, r/2 \}, \quad t \in \R \label{eq.orbitchg}
\end{equation}
for some $g \in \Rep(\min \{ \e_0/4, r/2 \})$ with $g(0) = 0$.
Now apply Lemma \ref{lem.orbitversion} to $x_1$ and $x_2$ respectively,
to have $h_1, h_2 \in \Rep(\e_0/4)$ with $h_1(0) = h_2(0) = 0$ such that
\begin{align}
\dist \bigl( \phi(t, x_0), \, \phi(h_1(t), x_1) \bigr) < \e_0/4, \quad t \geq 0,\\
\dist \bigl( \phi(t, x_0), \, \phi(h_2(t), x_2) \bigr) < \e_0/4, \quad t \geq 0.
\end{align}
Then these inequalities and \eqref{eq.orbitchg} show that for all $t \geq 0$,
\begin{align}
\dist \bigl( \xi(t + t_0), \,& \phi( h_1 \circ h_2^{-1} \circ g (t), x_1) \bigr)\\
&\leq \dist \bigl( \xi(t + t_0), \, \phi( g (t), x_2) \bigr)\\
&+\dist \bigl( \phi( g (t), x_2), \, \phi( h_2^{-1} \circ g (t), x_0) \bigr)\\
&+ \dist \bigl( \phi( h_2^{-1} \circ g (t), x_0), \, \phi( h_1 \circ h_2^{-1} \circ g (t), x_1) \bigr)\\
&< \frac{\e_0}{4} + \frac{\e_0}{4} + \frac{\e_0}{4} = \frac{3\e_0}{4}.
\end{align}
When $t \neq s$, we have
\begin{align}
\left| \frac{h_1 \circ h_2^{-1} \circ g (t) - h_1 \circ h_2^{-1} \circ g (s)}{t - s} \right|
= &\left| \frac{h_1 \circ h_2^{-1} \circ g (t) - h_1 \circ h_2^{-1} \circ g (s)}
{h_2^{-1} \circ g (t) - h_2^{-1} \circ g (s)} \right|\\
\cdot &\left| \frac{h_2^{-1} \circ g (t) - h_2^{-1} \circ g (s)}
{g (t) - g (s)} \right|\\
\cdot &\left| \frac{g (t) - g (s)}{t - s} \right|.
\end{align}
From this and $\e_0 < 4/5$ it follows that for all $t \in \R$,
\begin{align}
1 - \e_0
&<\left( 1 - \frac{\e_0}{4} \right) \left( 1 + \frac{\e_0}{4} \right)^{-1} \left( 1 - \frac{\e_0}{4} \right)\\
&< \left| \frac{h_2 \circ h_1^{-1} \circ g (t) - h_2 \circ h_1^{-1} \circ g (s)}{t - s} \right|\\
&< \left( 1 + \frac{\e_0}{4} \right) \left( 1 - \frac{\e_0}{4} \right)^{-1} \left( 1 + \frac{\e_0}{4} \right)
< 1 + \e_0.
\end{align}
Thus we have $h_1 \circ h_2^{-1} \circ g \in \Rep(\e_0)$, finishing the proof of Proposition \ref{pro.closed}.

\section{Neighborhoods of Singularities}\label{sec.neighborhoods_sing}
\begin{figure}[h]
\setlength{\unitlength}{20mm}
\begin{picture}(5, 5)
        \put(0.5, 0){\includegraphics[width=10cm,clip]{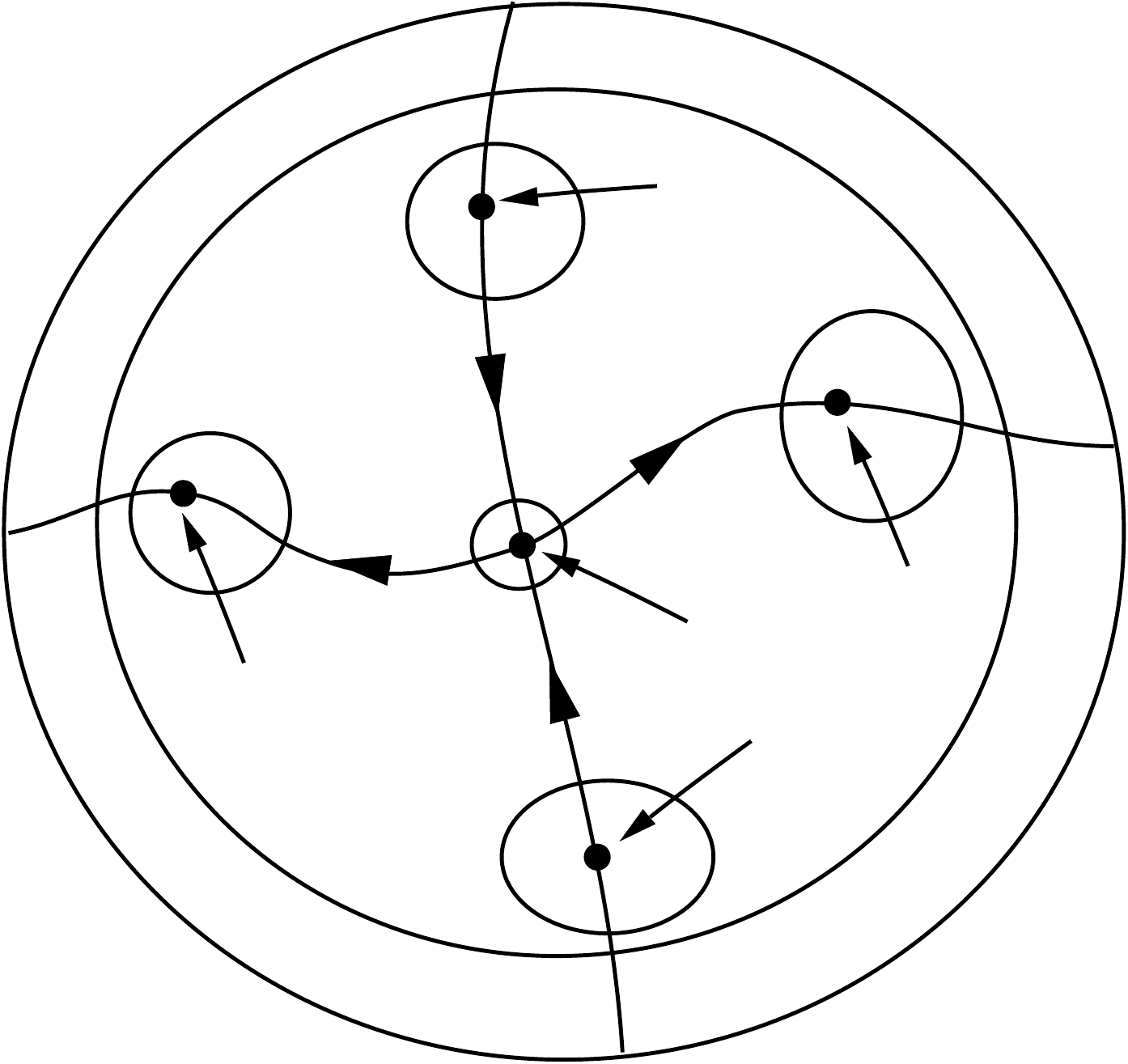}}
        \put(4.4, 4.4){\makebox(0, 0)[b]{\setlength{\fboxsep}{0.25mm}%
          \colorbox{white}{\normalsize $B(\e_0/2, p)$}}}
        \put(1.8, 4.0){\makebox(0, 0)[b]{\setlength{\fboxsep}{0.25mm}%
          \colorbox{white}{\normalsize $D$}}}
        \put(3.7, 1.8){\makebox(0, 0)[b]{\setlength{\fboxsep}{0.25mm}%
          \colorbox{white}{\normalsize $p$}}}
        \put(2.6,  2.4){\makebox(0, 0)[b]{\setlength{\fboxsep}{0.25mm}%
          \colorbox{white}{\normalsize $U_p$}}}
        \put(3.6, 3.8){\makebox(0, 0)[b]{\setlength{\fboxsep}{0.25mm}%
          \colorbox{white}{\normalsize $x_1$}}}
        \put(2.0, 3.4){\makebox(0, 0)[b]{\setlength{\fboxsep}{0.25mm}%
          \colorbox{white}{\normalsize $B(\tau_p, x_1)$}}}
        \put(4.0, 1.4){\makebox(0, 0)[b]{\setlength{\fboxsep}{0.25mm}%
          \colorbox{white}{\normalsize $x_2$}}}
        \put(2.4, 1.0){\makebox(0, 0)[b]{\setlength{\fboxsep}{0.25mm}%
          \colorbox{white}{\normalsize $B(\tau_p, x_2)$}}}
        \put(1.7, 1.6){\makebox(0, 0)[b]{\setlength{\fboxsep}{0.25mm}%
          \colorbox{white}{\normalsize $y_1$}}}
        \put(1.7, 2.7){\makebox(0, 0)[b]{\setlength{\fboxsep}{0.25mm}%
          \colorbox{white}{\normalsize $B(\tau_p, y_1)$}}}
        \put(4.6, 2.0){\makebox(0, 0)[b]{\setlength{\fboxsep}{0.25mm}%
          \colorbox{white}{\normalsize $y_2$}}}
        \put(4.0, 3.3){\makebox(0, 0)[b]{\setlength{\fboxsep}{0.25mm}%
          \colorbox{white}{\normalsize $B(\tau_p, y_2)$}}}
\end{picture}
\caption{}\label{figure.nbd42}
\end{figure}

In this section, we construct open sets near singularities in order to investigate the behavior of pseudotrajectories (see Figure \ref{figure.nbd42}).

Let $p \in \Sing_4(\phi) \cup \Sing_2(\phi)$.
Take a neighborhood $D \subset B(\e_0/2, p)$ of $p$ homeomorphic to a closed disk in $\R^2$.
Then $\partial D$ is a simple closed curve surrounding $p$.
If $p \in \Sing_4(\phi)$, take $x_1, x_2, y_1, y_2 \in D$ such that
the four hyperbolic sectors associated with $\partial D$ of $p$ is determined by basic orbits:
\[
\Orb^+(\phi(-1, x_i)) \text{ and } \Orb^-(\phi(1, y_j))
\]
with $i, j \in \{ 1, 2 \}$.
Similarly, if $p \in \Sing_2(\phi)$, take $x_1, y_1 \in D$ such that
the two hyperbolic sectors associated with $\partial D$ of $p$ is determined by basic orbits
\[
\Orb^+(\phi(-1, x_1)) \text{ and } \Orb^-(\phi(1, y_1)).
\]
Define
\[
G^s(p) = 
\begin{cases}
\{ x_1, x_2 \}, & p \in \Sing_4(\phi),\\
\{ x_1 \}, & p \in \Sing_2(\phi)
\end{cases}
\]
and
\[
G^u(p) = 
\begin{cases}
\{ y_1, y_2 \}, & p \in \Sing_4(\phi),\\
\{ y_1 \}, & p \in \Sing_2(\phi).
\end{cases}
\]
The following proposition claims that
the $\alpha$-limit set set of a point in $W^s(p) \setminus \{ p \}$ or
the $\omega$-limit set set of a point in $W^u(p) \setminus \{ p \}$
for $p \in \Sing_4(\phi) \cup \Sing_2(\phi)$ displays a simple dynamics.
More precisely, the following proposition shows that there is no saddle connection.
\begin{pro}\label{pro.connline}
Let $p \in \Sing_4(\phi) \cup \Sing_2(\phi)$.
If $z_0 \in W^s(p) \setminus \{ p \}$,
$($resp. $z_0 \in W^u(p) \setminus \{ p \}$$)$
then $\alpha(z_0)$ $($resp. $\omega(z_0)$$)$ is a backward asymptotically stable $($resp. asymptotically stable$)$ critical element.
\end{pro}
\begin{prf}
Suppose that $z_0 \in W^s(p) \setminus \{ p \}$.
Let $\Lambda_0$ be the critical element containing $\alpha(z_0)$.
Then we have that $\Lambda_0$ is not an asymptoticallly stable critical element.
Arguing by contradiction we assume that $\Lambda_0 \in \Sing_4(\phi) \cup \Sing_2(\phi)$ (see Theorem \ref{thm.classofsing} and Lemma \ref{lem.classclosedtraj}).
Let
\[
d_1 \in(0, \inf\{ L_0(\Lambda_0) - L_0(\Lambda) ; \Lambda \longrightharpoonup \Lambda_0 \}).
\]
Notice that $L_0(\Lambda_0) > L_0(\Lambda)$ for all $\Lambda$ with $\Lambda \longrightharpoonup \Lambda_0$ (see Lemma \ref{lem.ordering}).
Making $d_1$ smaller if necessary, we may take a neighborhood $N^u(p)$ of $G^u(p)$ such that $L_0(y) > L_0(p) + d_1/2$ for all $y \in N^u(p)$.
Then by Proposition \ref{pro.flowbox}, we may take an embedded closed interval $\Sigma_0$ satisfying the condition in Proposition \ref{pro.flowbox} with $x_0 = \phi(1, z_0)$.
By Lemma \ref{lem.hypsecnbd} we may assume that
\begin{equation}
\Orb^+(x)\cap N^u(p) \neq \emptyset \label{eq.closedprohyp}
\end{equation}
for all $x \in \Sigma_0 \setminus \{ \phi(1, z_0) \}$.
Take $\e_1 > 0$ so small that
\begin{equation}
\lvert L_0(x) - L_0(y) \rvert \leq d_1/2  \label{eq.100}
\end{equation}
for all $x, y \in M$ with $\dist(x, y) < \e_1$
and
\[
B(\e_1, z_0) \subset \{ \phi(t, x) ; -2 < t < 0, x \in \Sigma_0 \}.
\]
Since $\phi \in \OS(M)$, there exists $d > 0$ such that every $\xi \in \Pt(d)$ is $\e_1$-oriented shadowed.
Take $z_1 \in W^s(\Lambda_0) \setminus \Lambda_0$ and define
\[
\xi(t) = 
\begin{cases}
\phi(t, z_0), & t \geq -T_d,\\
\phi(t + 2T_0, z_1), & t < -T_d
\end{cases}
\]
for some $T_d > 0$.
Then $\xi \in \Ps(d)$ for sufficiently large $T_d$.
Thus there exist $h \in \Rep$ with $h(0) = 0$ and $w_0 \in M$ such that
\[
\dist(\xi(t), \phi(h(t), w_0)) < \e_1, \quad t \in \R.
\]
In particular, \eqref{eq.100} implies
\begin{equation}
\lvert L_0(\xi(t)) - L_0(\phi(h(t), y_0)) \rvert \leq d_1/2, \quad t \in \R.\label{eq.l0app}
\end{equation}
Since $w_0 \in B(\e_1, z_0)$, the choice of $\e_1$ implies that
$\phi(s, w_0) \in \Sigma_0$ for some $s \in [0, 2]$.
Now suppose that $\phi(s, w_0) \neq z_0$.
Then, by \eqref{eq.closedprohyp}, we may take $t' \geq 0$ satisfying
$\phi(h(t'), w_0) \in N^u(p)$,
which shows
\[
L_0(\phi(h(t'), w_0)) >  L_0(p) + d_1/2
\]
by the choice of $N^u(p)$.
On the other hand, it is obvious that
\[
L_0(p) = \sup \{ L_0(\phi(t, z_0)) ; t \in \R \} > L_0(\phi(t', z_0)) = L_0(\xi(t')),
\] 
contradicting inequality \eqref{eq.l0app} with $t = t'$.
This contradiction gives $\phi(s, w_0) = z_0$.
From this and \eqref{eq.l0app}, it follows that
\begin{eqnarray}
L_0(\Lambda_0) &=& 
\inf \{ L_0(\phi(t, z_0)) ; t \in \R \}\\ &=& 
\inf \{ L_0(\phi(t, w_0)) ; t \in \R \}\\ &\leq&
\inf \{ L_0(\phi(t, z_1)) ; t \in \R \} + d_1/2,
\end{eqnarray}
which contradicts $z_1 \in W^s(\Lambda_0) \setminus \Lambda_0$ and the choice of $d_1$.

The proof for the $\omega$-limit set is similar.
\qed
\end{prf}
The following proposition is the singularities version of Proposition \ref{pro.closed}.
However the proof is much simpler because $h \in \Rep(\e_0)$ can be the identity.
\begin{pro}\label{pro.sing}
Let $p \in \Sing_s(\phi)$ $($resp. $p \in \Sing_u(\phi)$$)$.
Then there exist a neighborhoood $U_p \subset B(\e_0, p)$ of $p$ and $d_p, r_p > 0$ 
such that if $\xi \in \Pt(d_p)$ satisfies $\xi(t_0) \in U_p$ for some $t_0 \in \R$ then the following properties hold:
\begin{enumerate}
  \item[(a)]
  $\xi(t) \in B(\e_0/2, p)$ for all $t \geq t_0$ $($resp. $t \leq t_0$$)$.
  \item[(b)]
  For any $x \in B(r_p, \xi(t_0))$,
  \[
  \dist \bigl( \xi(t), \, \phi(t - t_0, x) \bigr) < \e_0, \quad t \geq t_0
  \text{ $($resp. $t \leq t_0$$)$. }
  \]
\end{enumerate} 
\end{pro}
\begin{prf}
It is enough to consider $p \in \Sing_s(\phi)$ since the case for $p \in \Sing_u(\phi)$ is similar.
It follows from Lemma \ref{lem.genstable} that there exist $d > 0$ and a neighborhood $V_p$ of $p$ satisfying the following property:
\begin{itemize}
\item If $\xi \in \Ps(d)$ satisfies $\xi(0) \in V_p$, then $\xi(t) \in B(\e_0/2, p)$ for $t \geq 0$.
\end{itemize}
Take $r_p > 0$ so small that $B(2r_p, p) \subset V_p$.
Set $U_p = B(r_p, p)$ and let $d_p > 0$ be such that $\Pt(d_p) \subset \Ps(d)$.
By hypothesis, we have $\xi \in \Pt(d_p)$ with $\xi(t_0) \in U_p$ for some $t_0 \in \R$.
Since $U_p \subset V_p$, the choice of $V_p$ implies property $(a)$.
As for property $(b)$, by the choice of $r_p$, we have
$B(r_p, \xi(t_0)) \subset B(2r_p, p) \subset V_p$.
Thus, using the above property twice (notice that $t \mapsto \phi(t, x)$ is a $d$-pseudotrajectory), we have
\[
\dist \bigl( \phi(t - t_0, x), \, \xi(t) \bigr)
 < \dist \bigl( \phi(t - t_0, x), \, p \bigr)
 + \dist \bigl( p, \, \xi(t) \bigr) < \frac{\e_0}{2} + \frac{\e_0}{2} = \e_0,
\]
proving property $(b)$.
\qed
\end{prf}

For every $\Lambda \in \Closed_s(\phi) \cup \Closed_u(\phi) \cup \Sing_s(\phi) \cup \Sing_u(\phi)$,
fix a neighborhood $U_\Lambda$ of $\Lambda$ given by Proposition \ref{pro.closed} or Proposition \ref{pro.sing}.
In these propositions, we also gave constants $d_\Lambda, r_\Lambda > 0$.
Then let
\begin{equation}
d_0 = \min \{d_\Lambda ; \Lambda \in \Closed_s(\phi) \cup \Closed_u(\phi) \cup \Sing_s(\phi) \cup \Sing_u(\phi) \}\label{eq.defd_0}
\end{equation}
and
\begin{equation}
r_0 = \min \{r_\Lambda ; \Lambda \in \Closed_s(\phi) \cup \Closed_u(\phi) \cup \Sing_s(\phi) \cup \Sing_u(\phi) \}.\label{eq.defr_0}
\end{equation}

The following lemma claims that
any backward (resp. forward) pseudotrajectory with a point in the fundamental domain of $W^s(p)$ (resp. $W^u(p)$) with $p \in \Sing_4(\phi) \cup \Sing_2(\phi)$ can be $\e_0$-shadowed along the backward (resp. forward) pseudotrajectory.
\begin{lem}\label{lem.sing42center}
Let $p \in \Sing_4(\phi) \cup \Sing_2(\phi)$.
There exist $d_p, \tau_p > 0$ satisfying the following property:
\begin{itemize}
  \item
If $\xi \in \Pt(d_p)$ and $x_0 \in M$ satisfies $\xi(t_0), x_0 \in B(\tau_p, x)$ for some $t_0 \in \R$ and $x \in G^s(p)$ $($resp. $x \in G^u(p)$$)$,
then there exists $h \in \Rep(\e_0)$ with $h(0) = 0$ such that
\begin{equation}
\dist \bigl( \xi(t + t_0), \, \phi(h(t), x_0) \bigr) < \e_0,
\quad t \leq 0 \text{ $($resp. $t \geq 0$$)$ }.\label{eq.sing42center}
\end{equation}
\end{itemize}
\end{lem}
\begin{prf}
Let $x \in G^s(p)$.
Then Proposition \ref{pro.connline} shows that $\alpha(x) \in \Closed_u(\phi) \cup \Sing_u(\phi)$.
Take $T > 0$ satisfying $\phi(-T, x) \in U_{\alpha(x)}$.
Then choose $d_x \in (0, d_0)$ and $\tau_x > 0$ such that if $\xi \in \Pt(d_x)$ with $\xi(t_0) \in B(\tau_x, x)$ and $y \in B(\tau_x, x)$ then
\begin{equation}
\xi(t_0 - T) \in U_{\alpha(x)} \label{eq.disthaji}
\end{equation}
and
\begin{equation}
\dist \bigl( \xi(t_0 - t), \, \phi(-t, y) \bigr) < \min \{\e_0, r_{\alpha(x)} \}, \quad t \in [0, T]\label{eq.distint}.
\end{equation}

Suppose that $\xi \in \Pt(d_x)$ and $x_0$ satisfies $\xi(t_0), x_0 \in B(\tau_x, x)$.
By \eqref{eq.disthaji} and \eqref{eq.distint},
Propositions \ref{pro.closed} and \ref{pro.sing}
show that there exists $g \in \Rep(\e_0)$ with $g(0) = 0$ such that
\begin{equation}
\dist \bigl( \xi(t + t_0 - T), \, \phi(g(t), \phi(-T, x_0)) \bigr) < \e_0, \quad t \leq 0.\label{eq.defg}
\end{equation}
Define $h \in \Rep(\e_0)$ by
\[
h(t) =
\begin{cases}
g(t + T) - T, & t \leq -T,\\
t,  & t \geq -T.
\end{cases}
\]
Then, by \eqref{eq.distint},
\[
\dist \bigl( \xi(t + t_0), \, \phi(h(t), x_0) \bigr) < \e_0, \quad t \in [-T, 0].
\]
This inequality and \eqref{eq.defg} yield \eqref{eq.sing42center}.
Similarly, we can choose $d_y$ and $\tau_y$ for all $y \in G^s(p) \cup G^u(p)$.
Taking the minimum over the finite set  $G^s(p) \cup G^u(p)$, we complete the proof of the lemma.
\qed
\end{prf}
For $p \in \Sing_4(\phi) \cup \Sing_2(\phi)$ and $x \in G^s(p) \cup G^u(p)$,
fix an embedded closed interval
\begin{equation}
\Sigma_{x} \subset B(\tau_p, x)\label{eq.defofsigma}
\end{equation}
satisfying the condition of Proposition \ref{pro.flowbox} with $x_0 = x$.
When $x \in G^s(p)$ and $y \in G^u(p)$, there exists a point of $\Sigma_x$ which has the forward orbit intersecting $\Sigma_y$ at an arbitary large time.
\begin{lem}\label{lem.timenearsing}
Let $x_i, y_j \in D$ be the points given in the beginning of this section.
Then there exists $T_p > 0$ such that for all $i, j \in \{ 1, 2 \}$ and $T \geq T_p$,
\[
\phi(t, x) \in D, \quad t \in [0, T]
\]
and $\phi(T, x) \in \Sigma_{y_j}$
for some $x \in \Sigma_{x_i}$.
\end{lem}
\begin{prf}
Set $i, j \in \{ 1, 2 \}$ and let $S$ be the hyperbolic sector of $p$ determined by $\phi(-1, x_i)$ and $\phi(1, y_j)$.
By Lemma \ref{lem.hypsecnbd},
making $\Sigma_{x_i}$ smaller if necessary, we may assume that for any $x \in \Sigma_{x_i} \cap S$, there exists $t_x > 0$ such that
\[
\phi(t, x) \in S, \quad t \in [0, t_x]
\]
and $\phi(t_x, x) \in \Sigma_{y_j}$.
Since the map $\Sigma_{x_i} \cap S \ni x \mapsto t_x \in (0, \infty)$ is continuous, and $t_x \to \infty$ when $x \to x_i$, this lemma has been proved.
\qed
\end{prf}
The following lemma claims that some dynamical behavior around $p \in \Sing_4(\phi) \cup \Sing_2(\phi)$ is preserved for pseudotrajectories. 
\begin{lem}\label{lem.sing-c}
Let $p \in \Sing_4(\phi) \cup \Sing_2(\phi)$ and let $d_p, \tau_p > 0$ and
$T_p > 0$ be the constants given by Lemmas \ref{lem.sing42center} and \ref{lem.timenearsing}, respectively.
Then there exist a neighborhoood $U_p \subset D$ of $p$ and $d_p' \in (0, d_p)$ satisfying the following properties:
\begin{enumerate}
  \item[(a)] Assume that $\xi \in \Pt(d_p')$ satisfies $\xi(t_0) \in U_p$ and
  $\xi(t_0 + t_1) \notin B(\e_0/2, p)$ for some $t_1 > 0$.
  Then there exists $t_2 > T_p/2$ such that
  \begin{equation}
  \xi(t_0 + t_2) \in B(\tau_p, G^u(p)) \label{eq.singcfir}
  \end{equation}
  and
  \begin{equation}
  \xi(t_0 + t) \in B(\e_0/2, p), \quad t \in [0, t_2].\label{eq.singcsec}
  \end{equation}
  \item[(b)] Assume that $\xi \in \Pt(d_p')$ satisfies $\xi(t_0) \in U_p$ and
  $\xi(t_0 + t_1) \notin B(\e_0/2, p)$ for some $t_1 < 0$.
  Then there exists $t_2 < -T_p/2$ such that
  \[
  \xi(t_0 + t_2) \in B(\tau_p, G^s(p))
  \]
  and
  \[
  \xi(t_0 + t) \in B(\e_0/2, p), \quad t \in [t_2, 0].
  \]
\end{enumerate} 
\end{lem}
\begin{prf}
We prove this lemma for $p \in \Sing_4(\phi)$;
for the proof for $p \in \Sing_2(\phi)$ is essentially the same.
Let $\tau_p > 0$ be the contant given in Lemma \ref{lem.sing42center}.
Taking $\tau_p' \in (0, \tau_p)$, we may assume that
\[
\dist \bigl( p, \, G^s(p) \cup G^u(p) \bigr) > \tau_p'.
\]
If a neighborhood $U$ of $p$ and $d > 0$ are small enough,  
\begin{equation}
\xi(t + t_0) \in B(\e_0/2, p) \setminus
\overline{B(\tau_p', G^s(p) \cup G^u(p))}, \quad t \in [-T_p/2, T_p/2]\label{eq.defofU}
\end{equation}
for all $\xi \in \Pt(d)$ with $\xi(t_0) \in U$ for some $t_0 \in \R$.

Choose $r_1 \in (0, \tau_p'/2)$ 
with $D \subset B(\e_0/2 - r_1, p)$.
By Lemma \ref{lem.hypsecnbd}, there exists $r_2 \in (0, r_1)$ such that
if $x \in B(2r_2, p)$ and $\Orb^+(x) \not\subset D$
then
\begin{equation}
\phi(t, x) \in D, \quad t \in [0, t_x] \label{eq.intCfir}
\end{equation}
and
\begin{equation}
\phi(t_x, x) \in B(r_1, G^u(p)) \label{eq.intCsec}
\end{equation}
for some $t_x > 0$.
In addition, making $r_2$ smaller if necessary, we may assume that $B(r_2, p) \subset U$.
Take $d_p' > 0$ such that every $\xi \in \Pt(d_p')$ is $r_2$-oriented shadowed.
Set $U_p = B(r_2, p)$.

As in the hypothesis of property (a),
suppose that $\xi \in \Pt(d_p')$ satisfies $\xi(t_0) \in U_p = B(r_2, p)$ and $\xi(t_0 + t_1) \notin B(\e_0/2, p)$ for some $t_0 \in \R$ and $t_1 > 0$.
Then there exist $x_0 \in B(2r_2, p)$ and $h \in \Rep$ with $h(0) = 0$ such that 
\begin{equation}
\dist \bigl( \xi(t_0 + t), \, \phi(h(t), x_0) \bigr) < r_2, \quad t \in \R,
\label{eq.orientt0t}
\end{equation}
where $x_0 \in B(2r_2, p)$ follows from
\[
\dist ( x_0, p)
\leq \dist \bigl( x_0, \, \xi(t_0) \bigr)
+ \dist \bigl( \xi(t_0), \, p \bigr)
< r_2 + r_2.
\]
Moreover, inequality
\[
\dist \bigl( p, \, \phi(h(t_1), x_0) \bigr)
\geq \dist \bigl( p, \, \xi(t_0 + t_1) \bigr)
- \dist \bigl( \xi(t_0 + t_1), \, \phi(h(t_1), x_0) \bigr)
> \e_0/2 - r_2
\]
shows that $\phi(h(t_1), x_0) \notin D$.
Thus \eqref{eq.intCfir} and \eqref{eq.intCsec} with $x = x_0$ hold for some $t_{x_0} > 0$.
Then \eqref{eq.intCfir}, \eqref{eq.orientt0t} and the choice of $r_1$, $r_2$ with
$D \subset B(\e_0/2 - r_1, p) \subset B(\e_0/2 - r_2, p)$ imply that \eqref{eq.singcsec} for $t_2 = h^{-1}(t_{x_0})$.
Now \eqref{eq.orientt0t} and \eqref{eq.intCsec}
give
\begin{align}
\dist \bigl( \xi(t_0 + h^{-1}(t_{x_0})), \, G^u(p) \bigr)
&\leq \dist \bigl( \xi(t_0 + h^{-1}(t_{x_0})), \, \phi(t_{x_0}, x_0) \bigr)\\
&+ \dist \bigl( \phi(t_{x_0}, x_0), \, G^u(p) \bigr)\\
&< r_2 + r_1 < 2r_1 < \tau_p' < \tau_p.
\end{align}
This and \eqref{eq.defofU} show that $h^{-1}(t_{x_0}) > T_p/2$.
In addition, this implies
\eqref{eq.singcfir} with $t_2 = h^{-1}(t_{x_0})$.
The proof for property (b) is similar.
\qed
\end{prf}

\section{Proof of Proposition \ref{pro.finalpurpose}}\label{sec.prfof2dim_2}
In this section, let us finish the proof of Theorem \ref{thm.2dim} by showing Proposition \ref{pro.finalpurpose} given in Section \ref{sec.prolog_thm1_2}.
We need Lemma \ref{lem.shadowing} and Proposition \ref{pro.closed} and two propositions and three lemmas in Section \ref{sec.neighborhoods_sing}.
To apply Lemma \ref{lem.shadowing} later with $\e = \e_0$,
let $\e_0' > 0$ be the constant corresponding to $\e'$ in the proposition.
Let $d_0$ and $d_p'$ be constants given by \eqref{eq.defd_0} and Lemma \ref{lem.sing-c} respectively, and take
\[
d_1 < \min \{ d_0, d_p' ; p \in \Sing_4(\phi) \cup \Sing_2(\phi) \}
\]
so that every $\xi \in \Pt(d_1)$ is $\e_0'$-oriented shadowed.
Let
\[
U = \bigcup_{\Lambda \in \Critical(\phi)} U_\Lambda,
\]
where $U_\Lambda$ is the neighborhood of $\Lambda$ given in Propositions \ref{pro.closed}, 
\ref{pro.sing} and Lemma \ref{lem.sing-c}.

It is easy to see that the following claim follows from Lemma \ref{lem.birkhoffconst}.
\begin{claim}\label{claim.birkhoffconst}
There exist $R_0 > 0$ and $d_2 \in (0, d_1)$ such that
if $\xi \in \Pt(d_2)$ satisfies
\[
\xi(t) \notin U, \quad t \in (0, \ell),
\]
then $\ell < R_0$.
\end{claim}
In this claim, making $d_2 > 0$ smaller if necessary, we may assume that for all $\xi \in \Pt(d_2)$ and $t_0 \in \R$,
\begin{gather}
\dist \bigl( \xi(t_0 + t), \, \phi(t, \xi(t_0)) \bigr) < \min \{ \e_0, r_0 \}, \quad t \in [0, R_0], \label{eq0.03}
\end{gather}
where $\e_0$ and $r_0$ were given in \eqref{eq.defe_0} and \eqref{eq.defr_0}.

The hypothesis of the follwoing propsition is used for classsifying $\xi \in \Pt(d_2)$ into several cases.
\begin{pro}\label{pro.final}
Suppose that $\xi(t_0) \in U_p$ for some $p \in \Sing_4(\phi) \cup \Sing_2(\phi)$.
If
\[
\sup \{ t \in \R ; \xi(t) \in U_p \} = \infty
\text{ $($resp. $\inf \{ t \in \R ; \xi(t) \in U_p \} = -\infty$$)$},
\]
then $\xi(t) \in B(\e_0/2, p)$ for all $t \geq t_0$ $($resp. $t \leq t_0$$)$.
\end{pro}
\begin{prf}
Let us consider the former case.
Suppose to the contrary that there exists $t_1 > t_0$ with $\xi(t_1) \notin B(\e_0/2, p)$.
Then Lemma \ref{lem.sing-c} implies that
$\xi(t_0 + t_2) \in B(\tau_p, G^u(p))$ for some $t_2 > T_p/2$.
Then Proposition \ref{pro.connline} for asymptotically stable critical elements and Lemma \ref{lem.sing42center} for $t \geq 0$ show that
$\xi(t) \in B(3\e_0/2, \Lambda)$ for some asymptotically stable critical element $\Lambda$ and sufficiently large $t$, contradicting our hypothesis.
The proof for the latter case is similar.
\qed
\end{prf}
Let $\xi \in Pt(d_0)$ and define
\[
S(\xi) = \{ \Lambda \in \Critical(\phi) ; \xi(t) \in U_\Lambda \text{ for some $t \in \R$} \}.
\]


If $\xi(t_0) \in U_\Lambda$ for some $t_0 \in \R$ and some asymptotically stable critical element $\Lambda$, then by Propositions \ref{pro.sing} and \ref{pro.closed},
we have $\xi(t) \in B(\e_0, \Lambda)$ for all $t \geq t_0$.
Therefore $S(\xi)$ contains at most one asymptotically stable critical element.
Similarly, $S(\xi)$ contains at most one backward asymptotically stable critical element.
Thus it is enough to consider the following cases:
\begin{enumerate}
  \item[{\bf (C1)}] $S(\xi) = \{ \Lambda \}$ with $\Lambda \in \Sing_s(\phi) \cup \Closed_s(\phi)$.
  \item[{\bf (C2)}] $S(\xi) = \{ \Lambda \}$ with $\Lambda \in \Sing_u(\phi) \cup \Closed_u(\phi)$.
  \item[{\bf (C3)}] $S(\xi) = \{ \Lambda_1, \Lambda_2 \}$ with $\Lambda_1 \in \Sing_s(\phi) \cup \Closed_s(\phi)$ and $\Lambda_2 \in \Sing_u(\phi) \cup \Closed_u(\phi)$.
  \item[{\bf (C4)}] There exists $p \in (\Sing_4(\phi) \cup \Sing_2(\phi)) \cap S(\xi)$ such that
  \[
  \inf \{ t \in \R ; \xi(t) \in U_p \} = -\infty \text{ and }
  \sup \{ t \in \R ; \xi(t) \in U_p \} = \infty.
  \]
  \item[{\bf (C5)}] There exists $p \in (\Sing_4(\phi) \cup \Sing_2(\phi)) \cap S(\xi)$ such that
  \[
  \inf \{ t \in \R ; \xi(t) \in U_p \} > -\infty \text{ and }
  \sup \{ t \in \R ; \xi(t) \in U_p \} = \infty.
  \]
  \item[{\bf (C6)}] There exists $p \in (\Sing_4(\phi) \cup \Sing_2(\phi)) \cap S(\xi)$ such that
  \[
  \inf \{ t \in \R ; \xi(t) \in U_p \} = -\infty \text{ and }
  \sup \{ t \in \R ; \xi(t) \in U_p \} < \infty.
  \]
  \item[{\bf (C7)}] There exists $p \in (\Sing_4(\phi) \cup \Sing_2(\phi)) \cap S(\xi)$ such that
  \[
  \inf \{ t \in \R ; \xi(t) \in U_p \} > -\infty \text{ and }
  \sup \{ t \in \R ; \xi(t) \in U_p \} < \infty.
  \]
\end{enumerate}
For the proof of Theorem \ref{thm.2dim}, it suffices to show that any $\xi \in \Pt(d_2)$ is $\e_0$-standard shadowed in each case.
Now let us prove it for the seven cases.
\vspace{3mm}\\
{\bf Proofs for cases from (C1) to (C7)}.
\vspace{3mm}

{\bf (C1):}
From Claim \ref{claim.birkhoffconst}, it follows that $\{ t \in \R ; \xi(t) \in U_\Lambda \}$ does not have any lower bound.
Then, by Propositions \ref{pro.sing} and \ref{pro.closed},
\begin{equation}
\xi(t) \in B(\e_0/2, \Lambda), \quad t \in \R.\label{eq.caseC1}
\end{equation}
In fact, otherwise there is $t' \in \R$ such that $\xi(t') \notin B(\e_0/2, \Lambda)$,
but since there is no lower bounds, we can find some $t < t'$ such that $\xi(t) \in U_\Lambda$, which implies $\xi(t') \in B(\e_0/2, \Lambda)$ by these propositions, contradicting the above.
From \eqref{eq.caseC1} it follows that
,$\xi$ is $\e_0$-shadowed by $\phi(t, \Lambda)$ for all $t \in \R$ when $\Lambda \in \Sing_s(\phi)$.
On the other hand, when $\Lambda \in \Closed_s(\phi)$,
by the choice of $d_2$, there exist $x_0 \in M$ and $h \in \Rep$ with $h(0) = 0$ such that
\[
\dist \bigl( \xi(t), \, \phi(h(t), x_0) \bigr) < \e_0', \quad t \in \R.
\]
Then, by \eqref{eq.caseC1} and \eqref{eq.hyptildeK}, we can apply Lemma \ref{lem.shadowing} to
$[nT_0, (n + 1)T_0]$, $n \in \Z$, (thinking of $[0, T_0]$ in the proposition) in order to show that
$\xi$ is $\e_0$-standard shadowed.
\vspace{3mm}

{\bf (C2):}
The proof for this case is similar to that for case {\bf (C1)}.
\vspace{3mm}

{\bf (C3):}
By Claim \ref{claim.birkhoffconst}, we may take $t_1, t_2 \in \R$ such that
\[
\xi(t_1) \in U_{\Lambda_1}, \quad
\xi(t_2) \in U_{\Lambda_2}
\]
and $t_1 - t_2 \in (0, R_0)$.
By \eqref{eq0.03}, we have
\begin{equation}
\dist \bigl( \xi(t_2 + t), \, \phi(t, \xi(t_2)) \bigr) < \min \{ \e_0, r_0 \}, \quad t \in [0, R_0],\label{eq.finintquasi}
\end{equation}
implying
\[
\dist \bigl( \xi(t_1), \, \phi(t_1 - t_2, \xi(t_2)) \bigr) < r_0.
\]
Thus, Propositions \ref{pro.sing} and \ref{pro.closed} (b) show that
there exist $h_1, h_2 \in \Rep(\e_0)$ with $h_1(0) = h_2(0) = 0$
such that
\begin{equation}
\dist \bigl( \xi(t), \, \phi(h_1(t - t_1), \phi(t_1 - t_2, \xi(t_2))) \bigr) < \e_0, \quad t \geq t_1\label{eq.casec3fir}
\end{equation}
and
\begin{equation}
\dist \bigl( \xi(t), \, \phi(h_2(t - t_2), \xi(t_2)) \bigr) < \e_0, \quad t \leq t_2.\label{eq.casec3sec}
\end{equation}
Define $h \in \Rep(\e_0)$ with $h(t_2) = 0$ by
\[
h(t) = 
\begin{cases}
h_1(t - t_1) + t_1 - t_2, & t \geq t_1,\\
t - t_2, & t \in [t_2, t_1],\\
h_2(t - t_2), & t \leq t_2.\\
\end{cases}
\]
Then it follows from \eqref{eq.finintquasi}, \eqref{eq.casec3fir} and \eqref{eq.casec3sec} that
\[
\dist \bigl( \xi(t), \, \phi(h(t), \xi(t_2)) \bigr) < \e_0, \quad t \in \R.
\]
\vspace{3mm}

{\bf (C4):}
By Proposition \ref{pro.final}, $\xi(t) \in B(\e_0/2, p)$ for all $t \in \R$.
Therefore $\xi$ is $\e_0$-shadowed by $\phi(t, p)$ for $t \in \R$.
\vspace{3mm}

{\bf (C5):}
Assume that $\xi(t_0) \in U_p$.
By Claim \ref{claim.birkhoffconst} and $\inf \{ t \in \R ; \xi(t) \in U_p \} > -\infty$, we can find $t_1 < t_0$ with $\xi(t_1) \notin B(\e_0/2, p)$.
Then Lemma \ref {lem.sing-c} (b) implies that there exists $t_2 < -T_p/2$ such that
\[
\xi(t_0 + t_2) \in B(\tau_p, G^s(p))
\] 
and
\begin{equation}
\xi(t_0 + t) \in B(\e_0/2, p), \quad t \in [t_2, 0].\label{eq.t2t}
\end{equation}
Let $x_0 \in G^s(p)$ with $\xi(t_0 + t_2) \in B(\tau_p, x_0)$.
Then, by Lemma \ref{lem.sing42center}, there exists $h_0 \in \Rep(\e_0)$ with $h_0(0) = 0$ such that
\begin{equation}
\dist \bigl( \xi(t + t_0 + t_2), \, \phi(h_0(t), x_0) \bigr) < \e_0,
\quad t \leq 0.\label{eq.case5h0}
\end{equation}
Define $h \in \Rep(\e_0)$ with $h(0) = 0$ by
\[
h(t) =
\begin{cases}
h_0(t), & t \leq 0,\\
t, & t \geq 0.
\end{cases}
\]
Then Proposition \ref{pro.final}, \eqref{eq.t2t} and \eqref{eq.case5h0} imply
\[
\dist \bigl( \xi(t + t_0 + t_2), \, \phi(h(t), x_0) \bigr) < \e_0,
\quad t \in \R.
\]
\vspace{3mm}

{\bf (C6):}
The proof for this case is similar to that for case {\bf (C5)}.
\vspace{3mm}

{\bf (C7):}
Assume that $\xi(t_0) \in U_p$ for some $t_0 \in \R$.
By Proposition \ref{pro.final} and
the hypothesis,
there exist $t_1 < t_0$ and $s_1 > t_0$ with $\xi(t_1), \xi(s_1) \notin B(\e_0/2, p)$.
Then Lemma \ref{lem.sing-c} shows that there are $t_2 < -T_p/2$ and $s_2 > T_p/2$ such that
$\xi(t_0 + t_2) \in B(\tau_p, G^s(p))$, $\xi(t_0 + s_2) \in B(\tau_p, G^u(p))$ and
\begin{equation}
\xi(t_0 + t) \in B(\e_0/2, p), \quad t \in [t_2, s_2].\label{eq.nearsingpseudo}
\end{equation}
Let $x_0 \in G^s(p)$ and $y_0 \in G^u(p)$ be such that $\xi(t_0 + t_2) \in B(\tau_p, x_0)$ and
$\xi(t_0 + s_2) \in B(\tau_p, y_0)$, respectively.
By Lemma \ref{lem.timenearsing},
we can find $x_1$ in $\Sigma_{x_0}$ satisfying
\begin{equation}
\phi(t, x_1) \in B \subset B(\e_0/2, p), \quad t \in [0, s_2 - t_2]\label{eq.nearsingprbit}
\end{equation}
and $\phi(s_2 - t_2, x_1) \in \Sigma_{y_0}$,
where $\Sigma_{x_0}$ and $\Sigma_{y_0}$ are the embedded closed intervals given by \eqref{eq.defofsigma}.
Then it follows from Lemma \ref{lem.sing42center} that there exist $h_0 \in \Rep(\e_0)$ with $h_0(0) = 0$ such that
\begin{equation}
\dist \bigl( \xi(t + t_0 + t_2), \, \phi(h_0(t), x_1) \bigr) < \e_0,
\quad t \leq 0,\label{eq.case7h0}
\end{equation}
and
$g_0 \in \Rep(\e_0)$ with $g_0(0) = 0$ such that
\begin{equation}
\dist \bigl( \xi(t + t_0 + s_2), \, \phi(g_0(t), \phi(s_2 - t_2, x_1)) \bigr) < \e_0,
\quad t \geq 0. \label{eq.case7g0}
\end{equation}
Define $h \in \Rep(\e_0)$ with $h(0) = 0$ by
\[
h(t) =
\begin{cases}
h_0(t), & t \leq 0,\\
t, & t \in [0, s_2 - t_2],\\
g_0(t - s_2 + t_2) + s_2 - t_2, & t \geq s_2 - t_2.
\end{cases}
\]
Thus, by \eqref{eq.case7h0}, \eqref{eq.case7g0}, \eqref{eq.nearsingpseudo} and \eqref{eq.nearsingprbit}, we obtain
\[
\dist \bigl( \xi(t + t_0 + t_2), \, \phi(h(t), x_1) \bigr) < \e_0,
\quad t \in \R,
\]
completing the proof.

\section*{acknowledgements}
The author is grateful to my advisor S. Hayashi for his constructive suggestions and continuous support.

\bibliographystyle{amsplain}

\end{document}